\documentclass[preprint,12pt]{article}
\usepackage{epsfig,psfrag,verbatim,amssymb,amsfonts,amsthm}
\usepackage{graphicx}
\usepackage{subfigure,url,amsmath}
\usepackage[margin = 2.5cm, vmargin={2cm, 2cm},includefoot]{geometry}
\usepackage{hyperref}
\theoremstyle{plain}
\newtheorem{thm}{Theorem}[section]
\newtheorem{lem}{Lemma}[section]
\newtheorem{prop}{Proposition}[section]

\newtheorem{rem}{Remark}[section]
\newtheorem{defi}{Definition}[section]
\newtheorem{example}{Example}[section]

\newcommand{\BO} [1]{\mathbf{#1}}
\newcommand{\Bo} [1]{\boldsymbol{#1}}

\begin{document}
\title{Local Asymptotics of P-splines}

\author{Luo Xiao\thanks{Graduate student, Department of Statistical Science, Malott Hall,
Cornell University, New York 14853 (email: \texttt{lx42@cornell.edu}).},\
Yingxing Li\thanks{ Graduate student, Department of Statistical Science, Malott Hall,
Cornell University, New York 14853 (email: \texttt{yl377@cornell.edu}).},\
Tatiyana V. Apanasovich\thanks{Assistant Professor, Division of Biostatistics, Thomas Jefferson University,
Philadelphia, Pennsylvania 19107 (email: \texttt{Tatiyana.Apanasovich@jefferson.edu}).},\
and David Ruppert\thanks{Andrew Schultz, Jr., Professor of Engineering, School of
Operational Research and Information Engineering,  Comstock Hall, Cornell University,
New York 14853 (email: \texttt{dr24@cornell.edu}). } }

\date{June 7th, 2012}
\maketitle

\begin{abstract}
This report studies  local asymptotics of  P-splines with $p$th degree B-splines and a $m$th  order difference penalty.
Earlier work with $p$ and $m$ restricted is extended to the general case. Asymptotically, penalized splines are kernel estimators with equivalent kernels depending on $m$, but not on $p$. A central limit theorem provides simple expressions for the asymptotic mean and variance. Provided it is fast enough, the divergence rate of the number of knots does not affect the asymptotic distribution. The optimal convergence rate of the penalty parameter is given.
\vskip .2in

\end{abstract}
\newpage
\baselineskip 22pt

\numberwithin{equation}{section}

\numberwithin{equation}{section}

\section{Introduction}\label{sec:uni:introduction}
Suppose there is a univariate regression model
$$
y_i = \mu(x_i)+\epsilon_i, \quad i=1,\dots, n,
$$
where $\mu(x_i)$ and $\sigma^2(x_i)$ are the conditional expectation and variance of $y_i$ given $x_i$, respectively. For simplicity, we assume $x_i\in [0,1]$.

The regression function $\mu(x)$ can be modeled by $\sum_{k=1}^c\theta_k B_k(x)$ where $c=K+p$ and $\BO{B}(x) = \{B_1(x),\dots, B_c(x)\}^T$ is a  B-spline basis  of degree $p$ with knots $0=\kappa_0<\kappa_1<\dots<\kappa_{K} =1 $.  P-splines (Eilers and Marx, 1996)  find  $\hat{\Bo{\theta}}=(\hat{\theta}_1,\dots, \hat{\theta}_c)^T$  that minimizes 
\begin{equation}
\label{uni:est1}
\sum_{i=1}^n \left\{y_i - \sum_{k=1}^c \hat{\theta}_k B_k(x_i)\right\}^2+\lambda^{\ast} \sum_{k=m+1}^c\left \{\Delta^n\left(\hat{\theta}_k\right)\right\}^2, \quad \lambda \geq 0,
\end{equation}
where $\Delta$ is the difference operator, i.e., $\Delta(\theta_k)=\theta_k-\theta_{k-1}$ and $\Delta^m=\Delta(\Delta^{m-1})$, and $\lambda^{\ast}$ is the smoothing parameter.  Minimizing (\ref{uni:est1}) gives
\begin{align}
\label{uni:est2}
\left(\BO{B}^T\BO{B}/M+\lambda \mathbf{D}^T\mathbf{D}\right)\Bo{\hat{\theta}} = \BO{B}^T\BO{y},
\end{align}
where $M = n/K$, $\lambda = \lambda^{\ast}K/n$,  $\BO{y}=(y_1,\dots, y_n)^T$, $\BO{B}= \{\BO{B}(x_1)^T,\dots, \BO{B}(x_n)^T\}^T$ is an $n\times c$ matrix,  and $\mathbf{D}$ is the $m$th order differencing matrix of dimension  $(c-m)\times c$. For simplicity of notation, let 
\begin{equation}
\label{uni:Lambda}
\BO{\Lambda} = \BO{B}^T\BO{B}/M+\lambda \mathbf{D}^T\mathbf{D}
\end{equation} 
which is the smoother matrix for P-splines. Then the estimate  is given by 
\begin{align}
\label{uni:est3}
\hat{\mu}(x) =\BO{B}^T(x)\hat{\Bo{\theta}}= \BO{B}^T(x)\BO{\Lambda}^{-1}\BO{B}^T\BO{y}/M.
\end{align}

For simplicity, we  assume $x_1 = 1/(2n), x_2=3/(2n),  \dots, x_n=(2n-1)/(2n)$, i.e.,  the response is observed at equally spaced design points. We also assume $M$ is an integer to simplify some proofs. This assumption is for simplicity only and could be avoided. The case when the fixed design points are not equally spaced is considered in Section~\ref{sec:uni:irregular}.
\section{Review of Theoretical Study}\label{sec:uni:theory}
Penalized splines have been popular in recent years, as penalized splines use fewer  knots, thus need less computation than smoothing splines.  Ruppert {\it et al.}\ (2003) treat penalized splines extensively and also give numerous applications.

However, the theory of penalized splines has been remaining an interesting but challenging problem.  Opsomer and Hall (2005) first studied the asymptotic theory of penalized splines when $K$,  the number of knots, is infinite. Li and Ruppert (2008) derived the first  asymptotic distribution with low degree of splines and with low order of penalty. Wang {\it et al.}\ (2009) related penalized splines with some ordinary differential equations (ODEs), and by studying  Green's functions associated with those ODEs, they were able to derive the  asymptotic distribution of penalized splines. 

In contrast to Li and Ruppert (2008), Kauermann {\it et al.}\ (2009) considered the situation when $K$ increases at a moderate rate. Though they did not obtain an explicit expression for the asymptotic bias and variance, they generalized their results for non-normal responses. Claeskens {\it et al.}\ (2009) showed that depending on whether $K\rightarrow\infty$ increasing at a sufficiently fast or a sufficiently slow rate, the asymptotic distribution of penalized splines is either close to that of a smoothing spline or a regression spline. Correspondingly, they referred to these two cases as either  a large or small $K$ scenario. The large $K$ scenario is closest to current practice, as discussed, for example, in O'Sullivan (1986), Eilers and Marx (1996), and Ruppert {\it et al.}\ (2003), a relatively large number of knots is used and overfitting is controlled by a careful choice of smoothing parameter. 

One general approach to the theory of penalized splines is to use an equivalent kernel method,  which was first used by Silverman (1984) for studying the asymptotics of smoothing splines. The equivalent kernel method  was also useful in studying the asymptotics of P-splines (Li and Ruppert, 2008; Wang {\it et al.}, 2009). 

Independent from Wang {\it et al.}\ (2009), we extend Li and Ruppert's (2008) results  and provide an explicit expression on the asymptotic distribution of P-splines at an interior point. We  also derive the asymptotic distribution of P-splines near the boundary, acknowledging the existence of Wang {\it et al.}\ (2009). The conjecture,  that provided it is fast enough, the divergence rate of the number of knots does not affect the asymptotic distribution of penalized splines, is  confirmed in this paper.

The remainder of this chapter is organized as follows. In Section~\ref{sec:uni:main}, we summarize our main results. In Section~\ref{sec:uni:general}, we provide a general introduction of our method and present some technical results.  In Section~\ref{sec:uni:psplines_asymptotics}, We prove the main results in Section~\ref{sec:uni:main}. In Section~\ref{sec:uni:irregular}, we consider  irregularly spaced data. In Section~\ref{sec:uni:example}, we give an example illustrating the idea of binning data for irregularly space data. In Section~\ref{sec:uni:conclusion}, we conclude this chapter with some discussion. 

\section{Main Results}\label{sec:uni:main}
In this section, we summarize  the main results.  All derivations and proofs are given in Sections~\ref{sec:uni:general} and~\ref{sec:uni:psplines_asymptotics}. For notational convenience, $a\sim b$ implies $a/b$ converges to 1. We use the big ``O'' and small ``o" notation  that is with respect to $n$. Throughout this chapter, $ a = O(b)$ means $|a/b|$ converges to some finite nonnegative number as $n$ goes to infinity and $a = o(b)$ mean $|a/b|$ converges to 0. We also denote by $\mu^{(k)}(x)$ the $k$th derivative of the function $\mu(x)$. We need the following definition.
\begin{defi}
We define a kernel function
$$
H_m(x) = \frac{1}{2m}\sum_{\nu=1}^m \psi_{\nu} \exp\left(-\psi_{\nu}|x|\right),
$$
where $\psi_1,\cdots,\psi_m$ are the $m$ complex roots of $x^{2m}+(-1)^m=0$ such that all $\psi_{\nu} (1\leq \nu\leq m)$ have positive real parts.
\end{defi}

A kernel estimator with the kernel $H_m$ is of the form $(nh_n)^{-1}\sum_i y_i H_m\{h_n^{-1}(x-x_i)\}$, where $h_n$ is the bandwidth. As shown in Lemma~\ref{uni:lem9}, $H_m$ is of order $2m$ which determines the convergence rate the corresponding kernel estimator.  Proposition~\ref{uni:prop_est} shows that the P-spline estimator at an interior point is asymptotically equivalent to the above kernel estimator.

\begin{prop}
\label{uni:prop_est}
Assume the following conditions are satisfied.
\begin{enumerate}
\item There exists a constant $\delta>0$ such that
$\sup_{i}\textrm{\textnormal{E}}\left(|y_i|^{2+\delta}\right)<\infty$.  \item The regression function
$\mu(x)$ has a continuous $2m$th order derivative.
\item The variance function $\sigma^2(x)$ is
continuous. \item The random errors $\epsilon_{i}, 1\leq i\leq
n$,  are mutually independent. 
\item The covariates satisfy $x_i = (i-1/2)/n$, $1\leq i\leq n$.
\end{enumerate}
Let $\psi_0 = \min\{\text{Re}(\psi_1),\dots, \text{Re}(\psi_m)\}$, where $\text{Re}(\cdot)$ gives the real part of a complex number.  Let $h_n = \lambda^{1/(2m)}/K$.
Assume $h_n = o(1)$ and $(Kh_n)^{-1} = o(1)$. Let $\hat{\mu}(x)$ be the P-spline estimator using $m$th order difference penalty and $p$ degree B-splines with
equally spaced knots. Fix $x\in (0,1)$. Let $\mu^{\ast}(x) = (nh_n)^{-1}\sum_i y_i H_m\{h_n^{-1}(x-x_i)\}$. Then
\begin{equation*}
\begin{split}
\textrm{\textnormal{E}}\{\hat{\mu}(x)-\mu^{\ast}(x)\} &= O\left\{(Kh_n)^{-2}\right\},\\
\textrm{\textnormal{var}}\{\hat{\mu}(x)-\mu^{\ast}(x)\}& = o\left\{(nh_n)^{-1}\right\}.
\end{split}
\end{equation*}
\end{prop}

\begin{thm}\label{uni:thm1}
Use the same notation in Proposition~\ref{uni:prop_est} and assume all conditions and assumptions there are satisfied. Suppose that $K \sim  C n^{\tau}$ with $\tau>(m+1)/(4m+1)$, $h_n\sim hn^{-1/(4m+1)}$ for positive constants $C$ and $ h$ and $\lambda \sim (Kh_{n})^{2m}$. For any $x\in (0,1)$,
we have that
\[
n^{2m/(4m+1)}\left\{\hat{\mu}(x)-\mu(x)\right\}\Rightarrow N\left\{\tilde{\mu}(x),V(x) \right\}
\]
in distribution as $n\rightarrow\infty$,  where
\begin{align}
\tilde{\mu}(x) & = (-1)^{m+1}h^{2m}\mu^{(2m)}(x),\label{uni:tilde_mu(x)}\\
V(x) & =  \sigma^2(x)\int H_{m}^2(u)\mathrm{d}u.\label{uni:V(x)}
\end{align}
\end{thm}
\begin{rem}
Stone (1980) gave the optimal rates of convergence for nonparametric estimators.  For a univariate smooth function $\mu(x)$  with a continuous $2m$th  derivative, the corresponding optimal rate of convergence  for estimating $\mu(x)$ at any interior point is $n^{-2m/(4m+1)}$. Hence the P-spline estimator achieves the optimal rate of convergence.
\end{rem}

\begin{thm}\label{uni:thm2}
Assume conditions (1), (3), (4) and (5) in Proposition~\ref{uni:thm1} hold. Assume $\mu(x)$ has a continuous $m$th derivative over $[0,1]$. Suppose that $K \sim  C n^{\tau}$ with $\tau>(m+1)/(2m+1)$, $h_n\sim hn^{-1/(2m+1)}$ for positive constants $C$ and $ h$ and  $\lambda\sim (Kh_{n})^{2m}$. Let $\hat{\mu}(x)$ be the penalized estimator with $m$th order difference penalty and $p\geq 1$ degree B-splines with equally spaced knots. Assume $x\sim c_x h_n$ where $c_x$ is a constant. Then we have that
\[
n^{m/(2m+1)}\left\{\hat{\mu}(x)-\mu(x)\right\}\Rightarrow N\left\{\tilde{\mu}_0(x),V_0(x) \right\}
\]
in distribution as $n\rightarrow\infty$, where
\begin{align*}
\tilde{\mu}_0(x) & = (-1)^{m}h^{m}\mu^{(m)}(0)\int_{-\infty}^{c_x} u^{m}\left\{H_m(u)+H_{b,m}(c_x,c_x-u)\right\}\mathrm{d}u,\\
V_0(x) & =  \sigma^2(0)\int_{-\infty}^{c_x} \left\{H_{m}(u)+H_{b,m}(c_x, c_x-u)\right\}^2\mathrm{d}u.
\end{align*}
\end{thm}
Here $H_{b,m}$ is defined in~(\ref{uni:kernelb2}).
\begin{rem}
Theorems~\ref{uni:thm1} and~\ref{uni:thm2} show that the P-spline smoother has a slower rate of convergence at the boundary than in the interior.
\end{rem}

\section{Preliminary Derivation}\label{sec:uni:general}
We consider the large $K$ scenario (Claeskens {\it et al.}, 2009) and  assume $K$ and the smoothing parameter $\lambda$ increase with $n$ at certain rates specified later, respectively.  
 
The matrix $\BO{\Lambda}$ in~(\ref{uni:Lambda}) is a symmetric and banded matrix. For $q\leq k\leq  c- q$ with $q=\max(p,m)$, the $k$th column of $\BO{\Lambda}$ (denoted by $\BO{\Lambda}_k$) is
\[
(0,\dots, 0, \omega_q,\dots,\omega_1,\omega_0,\omega_1,\dots,\omega_q,0,\dots, 0)^T
\]
with the $k$th element being $\omega_0$. We need the following equation
\begin{equation}
\label{uni:P1}
\omega_q +\omega_{q-1}\rho+\cdots+\omega_1 \rho^{q-1}+\omega_0 \rho^q +\omega_1\rho^{q+1}+\cdots+\omega_q \rho^{2q}=0.
\end{equation}
Equation (\ref{uni:P1}) has a compact form
\begin{equation}
\label{uni:P2}
\lambda (-1)^m (1-\rho)^{2m}\rho^{q-m}+\rho^{q-p}P(\rho)=0,
\end{equation}
where 
\begin{equation}
\label{uni:P3}
P(x)= u_p +u_{p-1}x +\cdots +u_0 x^p +u_1 x^{p+1} +\cdots +u_p x^{2p}
\end{equation}
with the $k$th column of $\BO{B}^T\BO{B}$ being
\begin{equation}
\label{uni:P4}
(0,\dots,0,u_p,\dots, u_1,u_0,u_1,\dots, u_p,0,\dots, 0)^T.
\end{equation}
Let $\{\rho_{\nu}, \nu=1,\dots, q\}$ be the $q$ roots of (\ref{uni:P2}) such that when $\lambda$ is large, the real parts of  the first $m$  roots are all positive and less or equal than $1$ and moreover if $p>m$,  the other $q-m$ roots converge to zero. Define
\begin{equation}
\label{uni:S_k}
\mathbf{S}_k=\sum_{\nu=1}^q a_{\nu}\mathbf{T}_k(\rho_{\nu}),
\end{equation}
 where 
 \begin{equation}
 \label{uni:T_rho}
 \mathbf{T}_k (\rho)=(\rho^{k-1},\cdots,\rho,1,\rho,\cdots,\rho^{c-k})^T.
 \end{equation}
 For $1\leq \nu\leq q$ and  $2q\leq k\leq c-2q$, it can be shown that $\BO{T}_i(\rho_{\nu})$ is orthogonal to all  columns of $\BO{\Lambda}$ except the first $q$ columns, the last $q$ columns and the $j$th column with  $|k-j|< q$. The coefficient vector $\BO{a}=(a_1,\dots, a_q)^T$ can be chosen so that $\BO{S}_k$ is orthogonal to all  columns of $\BO{\Lambda}$ except the $k$th column,  the first $q$ columns and the last $q$ columns. It shall be shown later in this section that $\BO{a}$ does not depend on $k$. Specifically, we find a unique $\BO{a}$ such that
\begin{equation}
\label{uni:constraint1}
\mathbf{S}_k^T \BO{\Lambda}_k=1\quad \mbox{ and} \quad \mathbf{S}_k^T \BO{\Lambda}_j=0, \quad  0< |k-j|\leq q-1,
\end{equation}
where $\BO{\Lambda}_k$ is the $k$th column of $\BO{\Lambda}$ as before. 

Fix $x\in (0,1)$. By~(\ref{uni:est3}), we need only to consider non-zero $B_k(x)$. Hence we assume $k\in (Kx-p-1, Kx+p+1)$. By~(\ref{uni:constraint1}) and the definition of $\mathbf{S}_k$, there exists a constant $C>0$ such that,
\begin{equation}
\label{uni:e_k}
\BO{S}_k^T\BO{\Lambda}_j=   O\left[\exp\left\{-C \lambda^{-1/(2m)}K\min(x,1-x)\right\}\right],  \quad1\leq j\leq q, \mbox{and}\, c-q\leq j\leq c.
\end{equation} 
Let $\BO{e}_k$ be a vector of length $c$ with the $k$th entry $1$  and other elements 0. Define $\tilde \theta_k = (\mathbf{S}_k^T\BO{\Lambda})\hat{\boldsymbol{\theta}}$. Equation~(\ref{uni:est2}) implies
$\tilde{\theta}_k = \BO{S}_k^T\BO{B}^T\BO{y}$.
By~(\ref{uni:constraint1}), ~(\ref{uni:e_k}) and Lemma~\ref{uni:lem00}, 
$
\tilde\theta_k-\hat\theta_k = (\mathbf{S}_k^T\BO{\Lambda}-\BO{e}_k^T)\hat{\boldsymbol{\theta}}= \sum_{i=1}^n \tilde{b}_{i,k} y_i,
$
where $\tilde{b}_{i,k} =  O\left[\exp\left\{-C \lambda^{-1/(2m)}K\min(x,1-x)\right\}\right]$. Let $S_{k,r}$ be the $k$th element of $\BO{S}_k$. By (\ref{uni:est3}),
\begin{align}
\hat{\mu}(x) &= \sum_{k=1}^c B_k(x) \BO{S}_k^T\BO{B}^T\BO{y} +\sum_{k=1}^c B_k(x) (\tilde\theta_k -\hat\theta_k)\nonumber\\
&=\sum_{k=1}^c \left[B_k(x)\left\{\sum_{r=1}^c S_{k,r}\sum_{i=1}^n B_r(x_i)y_i\right\}\right]+\sum_{ |k-Kx|\leq p} B_k(x)\left(\sum_{i=1}^n \tilde{b}_{i,k}y_i\right) \nonumber\\
\label{uni:est5}
&=\sum_{i=1}^n y_i \left \{\sum_{k,r} B_k(x)B_r(x_i)S_{k,r}+b_i(x)\right\},
\end{align}
where $b_i(x) = \sum_{|k-Kx|\leq p} B_k(x)\tilde{b}_{i,k} = O\left[\exp\left\{-C \lambda^{-1/(2m)}K\min(x,1-x)\right\}\right]$.
 We assume appropriate regularity conditions on the data $\mathbf{y}$ so that interchanging sums in~(\ref{uni:est5}) is valid. Note that $\sum_{k,r} B_k(x)B_r(x_i)S_{k,r}+b_i(x)$ in (\ref{uni:est5}) is the weight of the $i$th observation for estimating $\hat{\mu}(x)$. 

For the boundary case, assume  $x$ goes to 0 at a rate of $\lambda^{1/(2m)}/K$, i.e., $x\sim c_x \lambda^{1/(2m)}/K$, where $c_x$ is  a constant. We assume that $\lambda^{1/(2m)}/K$ converges to 0.  Assume $k\in (Kx-p-1, Kx+p+1)$, then $\mathbf{S}_k$ is orthogonal to all columns of $\BO{\Lambda}$ except the $k$th, the first $q$ and the last $q$ columns. Furthermore, $\BO{T}_1(\rho)$ defined in~(\ref{uni:T_rho})  can be shown orthogonal to all columns of $\BO{\Lambda}$ except the first $q$ and  the last $q$ columns. 
Define $\BO{R}_k = \sum_{\nu=1}^q \tilde{a}_{k, \nu} \BO{T}_1(\rho_{\nu})$. Then $\BO{S}_k+\BO{R}_k$ is orthogonal to all columns of $\BO{\Lambda}$ except the $k$th, the first $q$ and  the last $q$ columns for arbitrary coefficient vector $\tilde{\BO{a}}_k=\{\tilde{a}_{k,1},\dots, \tilde{a}_{k, q}\}^T$. We find the coefficient vector $\tilde{\BO{a}}_k$  so that $\BO{S}_k+\BO{R}_k$ is orthogonal to all  columns of $\BO{\Lambda}$ except the $k$th and the last $q$ columns. Specifically, we find $\tilde{\BO{a}}$ such that
\begin{equation}
\label{uni:constraint2}
\left(\BO{S}_k+\BO{R}_k\right)^T \BO{\Lambda}_k=1\quad \mbox{ and} \quad \left(\BO{S}_k+\BO{R}_k\right)^T \BO{\Lambda}_j=0, \quad  0< j\leq  c-q.
\end{equation}
Then there exists a constant $C_0>0$ such that  for  $c-q\leq j\leq c$,  $(\BO{S}_k+\BO{R}_k)^T\BO{\Lambda}_j=O\left[\exp\left\{-C_0 \lambda^{-1/(2m)}K\right\}\right]$. We can derive that, similar to (\ref{uni:est5}), 
\begin{equation}
\label{uni:est6}
\hat{\mu}(x)=\sum_{i=1}^n y_i \left \{\sum_{k,r} B_k(x)B_r(x_i)(S_{k,r}+R_{k,r})+b_{i,0}(x)\right\},
\end{equation}
where $R_{k,r}$ is the $r$th element of $\BO{R}_k$ with $R_{k,r}= \sum_{\nu=1}^q \tilde{a}_{k, \nu}\rho_{\nu}^{r-1}$, and $b_{i,0}(x) = O\left[\exp\left\{-C_0 \lambda^{-1/(2m)}K\right\}\right]$.

In the next subsections,  we shall  derive  the coefficients $\rho_{\nu}, a_{\nu}$  and $\tilde{a}_{k, \nu}$. 
\subsection{Derivation of $\rho_{\nu}$}
\subsubsection{The case $p\leq m$}
In this case $q=m$. Equation (\ref{uni:P2}) becomes
\begin{equation}
\lambda(-1)^m (1-\rho)^{2m}+\rho^{m-p}P(\rho)=0
\label{uni:p<m}
\end{equation}
and $\rho_1,\dots,\rho_m$ are the  $m$ complex roots of (\ref{uni:p<m}) such that the real part of $\rho_{\nu}$ is positive and  less or equal than $1$.   Proposition~\ref{uni:prop_rho} below shows that $\rho_{\nu}$ exists and has  an explicit form. 
\begin{prop}
\label{uni:prop_rho}
As $\lambda\rightarrow \infty$, the roots of equation~(\ref{uni:p<m}) take the following forms
\begin{equation}
\label{uni:rho}
\rho_{\nu} = 1-\psi_{\nu} \lambda^{-1/(2m)} +1/2 \psi_{\nu}^2 \lambda^{-1/m} +O\left\{\lambda^{-3/(2m)}\right\}, \quad 1\leq\nu\leq 2m,
\end{equation}
where $\psi_1,\cdots,\psi_{2m}$ are the roots of $x^{2m}+(-1)^m=0$.
\end{prop}

\begin{rem}
To be consistent with the definition in Section~\ref{sec:uni:main}, we assume for the first $m$ roots, $\psi_{\nu}$ have positive real parts and for the last $m$ roots, $\psi_{\nu}$ have negative real parts. The real parts of $\rho_1,\dots, \rho_m$ are hence positive and equal or less than 1.
\end{rem}
{\it Proof of Proposition~\ref{uni:prop_rho}:} The existence of $2m$ roots for equation~(\ref{uni:p<m}) is obvious from complex analysis. Suppose $1-\delta_1$ is a root of equation~(\ref{uni:p<m}). Then 
$$
G_{1,\lambda}(\delta_1)=\lambda(-1)^m\delta_1^{2m}+(1-\delta_1)^{m-p}P(1-\delta_1)=0.
$$
Because the leading coefficient for the polynomial $G_{1,\lambda}(\delta_1)$ is $\lambda(-1)^m$ (or $\lambda(-1)^m+\omega_0$ if $m=p$), it is easy to see that $\delta_1$ is uniformly bounded as $\lambda\rightarrow\infty$. Hence $(1-\delta_1)^{m-p}P(1-\delta_1)$ is uniformly bounded, which implies $\lambda(-1)^m\delta_1^{2m}$ is uniformly bounded. It follows that $\lim_{\lambda\rightarrow\infty}\delta_1=0$. Then 
$$
 \lim_{\lambda\rightarrow\infty} G_{1,\lambda}(\delta_1)=\lim_{\lambda\rightarrow\infty}\lambda(-1)^m\delta_1^{2m} +1=0,
$$
which implies
\begin{equation}
\label{delta2}
\delta_1 = \psi_{\nu}\lambda^{-1/(2m)} (1+\delta_2),
\end{equation}
where $\psi_{\nu}$ is a root of $x^{2m}+(-1)^m=0$ for some $\nu$ and  $\lim_{\lambda\rightarrow\infty}\delta_2=0$. Substituting~(\ref{delta2}) into $G_{1,\lambda}$ (denoted by $G_{2,\lambda}(\delta_2)$) gives
\begin{equation}
\label{delta3_0}
0= G_{2,\lambda}(\delta_2)=-(1+\delta_2)^{2m}+\left\{1-\psi_{\nu}\lambda^{-1/(2m)}(1+\delta_2)\right\}^{m-p}P\left\{1-\psi_{\nu}\lambda^{-1/(2m)}(1+\delta_2)\right\}.
\end{equation}
It is easy to show that
\begin{align}
\label{delta3_1}
\left\{1-\psi_{\nu}\lambda^{-1/(2m)}(1+\delta_2)\right\}^{m-p}=&1-(m-p)\psi_{\nu}\lambda^{-1/(2m)}+o\left\{\lambda^{-1/(2m)}\right\},\\
\label{delta3_2}
P\left\{1-\psi_{\nu}\lambda^{-1/(2m)}(1+\delta_2)\right\} = &P(1)-P^{\prime}(1)\psi_{\nu}\lambda^{-1/(2m)}+o\left\{\lambda^{-1/(2m)}\right\}.
\end{align}
Equalities (\ref{delta3_0})--(\ref{delta3_2}), as well as Lemma~\ref{uni:lem0}, imply
$$
\delta_2 = \frac{p-m-P^{\prime}(1)}{2m}\psi_{\nu}\lambda^{-1/(2m)}(1+\delta_3) = - \frac{1}{2}\psi_{\nu}\lambda^{-1/(2m)}(1+\delta_3),
$$
where $\lim_{\lambda\rightarrow\infty}\delta_3=0$. By similar analysis, we can show that $\delta_3 = O\left\{\lambda^{-3/(2m)}\right\}$. Hence a root of equation~(\ref{uni:p<m}) takes the form
$$
1-\psi_{\nu} \lambda^{-1/(2m)} +1/2 \psi_{\nu}^2 \lambda^{-1/m} +O\{\lambda^{-3/(2m)}\}, \quad\text{for some }\,\nu.
$$
Thus,  equation~(\ref{uni:p<m}) has $2m$ roots that take the above form and each root has a $\psi_{\nu}$ that is a root of~(\ref{uni:rho}). 
\subsubsection{The case $p>m$}
When $p>m$, equation (\ref{uni:P2}) becomes
\begin{equation}
\label{uni:p>m}
\lambda(-1)^m (1-\rho)^{2m}\rho^{p-m} + P(\rho)=0.
\end{equation}
Similar to Proposition~\ref{uni:prop_rho}, we have the following
\begin{prop}
\label{uni:prop_rho_1}
As $\lambda\rightarrow \infty$, $2m$ roots of equation~(\ref{uni:p>m})  take the  forms in~(\ref{uni:rho}), and additionally, $p-m$ roots of equation~(\ref{uni:p>m}) take the following forms
\begin{equation}
\label{uni:rho_1}
\rho_{\nu}= \left\{\frac{\omega_q}{\lambda}\right\}^{\frac{1}{p-m}} \psi_{\nu} +O(\lambda^{-\frac{2}{p-m}}), \quad m+1\leq \nu\leq p,
\end{equation}
where $\psi_{m+1},\cdots,\psi_{p}$ are the roots of $x^{p-m}+(-1)^m=0$.
\end{prop}

{\it Proof of Proposition~\ref{uni:prop_rho_1}:}  Assume $\delta_0$ is a root of equation~(\ref{uni:rho_1}). Consider the case $\limsup_{\lambda\rightarrow\infty}\delta_0\neq 0$ and is bounded. Then a similar proof as that of Proposition~\ref{uni:prop_rho} gives $2m$ roots taking the forms in~(\ref{uni:rho}). Now consider the case $\limsup_{\lambda\rightarrow\infty}\delta_0= 0$.  $P(\delta_0)$ converges to $\omega_q$ as $\lambda\rightarrow\infty$, which implies $\lambda(-1)^m\delta_0^{p-m}$ converges to $-\omega_q$. It follows that $\delta_0 = \psi_{\nu}(\omega_q/\lambda)^{1/(p-m)}(1+\delta_1)$,  where $\psi_{\nu}$ is a root of $x^{p-m}+(-1)^m=0$ for some $\nu$ and $\lim_{\lambda\rightarrow\infty}\delta_1=0$.  Similar derivation as in the proof of Proposition~\ref{uni:prop_rho}  gives~(\ref{uni:rho_1}). To complete the proof, notice that for the case $\limsup_{\lambda\rightarrow\infty}\delta_0=\infty$, we can derive the rest $p-m$ unbounded roots of equation~(\ref{uni:p>m}).

\subsection{Derivation of $a_{\nu}$}
In this subsection, we shall establish the following
\begin{prop}
\label{uni:prop_anu}
Assume $q<k<c-q$ and $x\in(0,1)$. As $\lambda\rightarrow \infty$, the vector $\mathbf{a}$ satisfying the constraints in~(\ref{uni:constraint1}) is unique, i.e., does not depend on $k$, and has the following form
\begin{equation}
\label{uni:a2}
a_{\nu} = \frac{\psi_{\nu}}{2m}\lambda^{-1/(2m)}\left\{1+O(\lambda^{-1/m})\right\},  \quad 1\leq \nu \leq m,\end{equation}
and if $p>m$,
$$
a_{\nu} = O\left\{\lambda^{p/(m-p)}\right\}, \quad \nu=m+1,\dots, p.
$$
\end{prop}
\begin{rem}
Because the proof is lengthy, we shall sketch the proof within the context in the remainder of this subsection.
\end{rem}
 For $1\leq \nu \leq q$, define $s_j(\rho_{\nu})=\mathbf{T}_{k}^T(\rho_{\nu})\BO{\Lambda}_{i-q+j}$ for $1\leq j\leq q$. Then $s_j(\rho_{\nu})=\sum_{l=0}^{j-1} \omega_{q-l}(\rho_{\nu}^{j-l}-\rho_{\nu}^{l-j})$.  Constraints in~(\ref{uni:constraint1}) give a system of linear equations 
\begin{equation*}
\left(
\begin{array}{cccc}
  s_1(\rho_1) &\cdots & s_1(\rho_q)  \\
\vdots&\ddots&\vdots\\
s_{q-1}(\rho_1)&\cdots &s_{q-1}(\rho_q)\\
s_q(\rho_1)  & \cdots&  s_q(\rho_q)
\end{array}
\right)
\left(\begin{array}{c}
a_1\\
\vdots\\
a_{q-1}\\
a_q
\end{array}\right)=
\left(\begin{array}{c}
0\\
\vdots\\
0\\
1
\end{array}\right).
\end{equation*}
As shall be shown soon,  $a_{\nu}$'s exist and are unique.  Making use of the structure of $s_j(\rho_{\nu})$ and doing row transforms on the above linear equations, we have 
\begin{equation*}
\left(
\begin{array}{ccc}
 \omega_q(\rho_1-\rho_1^{-1}) &\cdots & \omega_q(\rho_q-\rho_q^{-1})  \\
 \vdots&\ddots&\vdots\\
\omega_q(\rho_1^{q-1}-\rho_1^{1-q}) &\cdots & \omega_q(\rho_q^{q-1}-\rho_q^{1-q})  \\
\omega_q(\rho_1^q-\rho_1^{-q}) &\cdots & \omega_q(\rho_q^q-\rho_q^{-q})  \\
\end{array}
\right)
\left(\begin{array}{c}
a_1\\
\vdots\\
a_{q-1}\\
a_q
\end{array}\right)=
\left(\begin{array}{c}
0\\
\vdots\\
0\\
1
\end{array}\right).
\end{equation*}
Further row transforms on the above equations give
\begin{equation*}
\left(
\begin{array}{ccc}
1&\cdots &1  \\
 \vdots&\ddots&\vdots\\
(\rho_1+\rho_1^{-1}-2)^{q-2}&\cdots &(\rho_q+\rho_q^{-1}-2)^{q-2}\\
(\rho_1+\rho_1^{-1}-2)^{q-1}&\cdots &(\rho_q+\rho_q^{-1}-2)^{q-1}\\
\end{array}
\right)
\left(\begin{array}{c}
a_1(\rho_1-\rho_1^{-1})\\
\vdots\\
a_{q-1}(\rho_{q-1}-\rho_{q-1}^{-1})\\
a_q(\rho_q-\rho_{q}^{-1})
\end{array}\right)=
\left(\begin{array}{c}
0\\
\vdots\\
0\\
\omega_q^{-1}
\end{array}\right).
\end{equation*}
In the above equations, the matrix before the column of coefficients is  a $q\times q$ Vandermonde matrix. Making use of the determinant property of Vandermonde matrix, the solution to the above linear equations exists and is unique because $\rho_{\nu}+\rho_{\nu}^{-1}-2,1\leq \nu\leq q$ are all different. Furthermore, it is apparent that the solution to the above equations does not depend on $k$, hence  $\BO{a}$ is the same for all $k$ such that $q\leq k\leq c-q$.  By Cramer's rule in solving  linear equations, we  obtain for $1\leq \nu\leq q$

\begin{equation}
\begin{split}
a_{\nu}\omega_q (\rho_{\nu}-1/\rho_{\nu})=&\frac{(-1)^{m+\nu}\prod_{1\leq i<j\leq q, j\neq \nu,i\neq \nu }(\rho_j+\rho_j^{-1}-\rho_i-\rho_i^{-1})}
{\prod_{1\leq i<j\leq q} (\rho_j+\rho_j^{-1}-\rho_i-\rho_i^{-1})}\\
=&\frac{(-1)^{q+\nu}(-1)^{q-\nu}}{\prod_{1\leq j\neq \nu \leq q}(\rho_{\nu}+\rho_{\nu}^{-1}-\rho_j-\rho_j^{-1})}\\
=&\frac{1}{\prod_{1\leq j\neq \nu\leq q}(\rho_{\nu}+\rho_{\nu}^{-1}-\rho_j-\rho_j^{-1})}.
\end{split}
\label{uni:eq2}
\end{equation}
Hence 
\begin{equation}
\label{uni:a1}
a_{\nu}^{-1}=\omega_q (\rho_{\nu}-\rho_{\nu}^{-1})\prod_{1\leq j\neq \nu \leq q}(\rho_{\nu}+\rho_{\nu}^{-1}-\rho_j-\rho_j^{-1}).
\end{equation}
\subsubsection{The case $p\leq m$}
By (\ref{uni:rho}), for $1\leq \nu \leq m$,
\begin{equation*}
\rho_{\nu}-\rho_{\nu}^{-1} = -2\psi_{\nu} \lambda^{-1/2m} +O(\lambda^{-3/2m}),
\end{equation*}
and 
\begin{equation*}
\rho_{\nu}+\rho_{\nu}^{-1}-2 = \psi_{\nu}^2\lambda^{-1/m} +O(\lambda^{-2/m}).
\end{equation*}
It follows that  for $1\leq j\neq \nu \leq m$,
\begin{equation}
\label{uni:eq3}
\rho_{\nu}+\rho_{\nu}^{-1}-\rho_j-\rho_j^{-1} = (\psi_{\nu}^2-\psi_j^2) \lambda^{-1/m}+O(\lambda^{-2/m}).
\end{equation}
Then
\begin{equation}
\begin{split}
\prod_{j\neq \nu}(\rho_{\nu}+\rho_{\nu}^{-1}-\rho_j-\rho_j^{-1})=&\lambda^{-1+1/m} \prod_{j\neq \nu}\left\{(\psi_{\nu}^2-\psi_j^2) +O(\lambda^{-1/m})\right\}\\
=&\lambda^{-1+1/m} \left\{\prod_{j\neq \nu}(\psi_{\nu}^2-\psi_j^2)+O(\lambda^{-1/m})\right\}.
\end{split}
\label{uni:eq1}
\end{equation}
By Lemma~\ref{uni:lemroot}, equality (\ref{uni:eq1}) can be simplified 
\begin{equation}
\label{uni:eq1_2}
\prod_{j\neq \nu}(\rho_{\nu}+\rho_{\nu}^{-1}-\rho_j-\rho_j^{-1})\\
=(-1)^{m+1}m  \psi_{\nu}^{-2}\lambda^{-1+1/m}\{1+O(\lambda^{-1/m})\}.
\end{equation}
In light of~(\ref{uni:eq2}) and~(\ref{uni:eq1_2}),
\begin{equation*}
\left\{a_{\nu}\omega_q (\rho_{\nu}-\rho_{\nu}^{-1})\right\}^{-1}=(-1)^{m+1}m^{-1}\psi_{\nu}^2 \lambda^{1-1/m}\{1 +O(\lambda^{-1/m})\}.
\end{equation*}
Note that  for $p\leq m$, $\omega_q=\omega_m = (-1)^m\lambda +$  a constant, where the constant is the coefficient  of $\rho^m$ in the polynomial  $P(\rho)$. Hence $(-1)^m \lambda^{-1}\omega_q =1+O(\lambda^{-1})$.  It follows that
\begin{equation*}
\begin{split}
a_{\nu}^{-1}=& \omega_q (\rho_{\nu}-1/\rho_{\nu})\prod_{j\neq \nu}(\rho_{\nu}+1/\rho_{\nu}-\rho_j-1/\rho_j)\\
=&-\omega_q \left\{ 2\psi_{\nu} \lambda^{-1/(2m)}+O(\lambda^{-3/(2m)})\right\} (-1)^{m+1}m \psi_{\nu}^{-2}  \lambda^{-1+1/m} \left\{1+O(\lambda^{-1/m})\right\}\\
 =&2m (-1)^m\lambda^{-1+1/(2m)}\omega_q \psi_{\nu}^{-1} \left\{1 +O(\lambda^{-1/m}) \right\}\\
 =&2m\lambda^{1/(2m)}\psi_{\nu}^{-1}\left\{1 +O(\lambda^{-1/m})\right\}.
\end{split}
\end{equation*}
The above derivation establishes~(\ref{uni:a2}).

\subsubsection{The case $p>m$}

To derive $a_{\nu}$, we need to study (\ref{uni:a1}) again. For  the term $\rho_{\nu}+\rho_{\nu}^{-1}-\rho_j-\rho_j^{-1}$ in (\ref{uni:a1}), there are two new cases besides~(\ref{uni:eq3}),
\begin{equation*}
\rho_{\nu}+\rho_{\nu}^{-1}-\rho_j-\rho_j^{-1} =\begin{cases}
 -\psi_j^{-1}(\lambda/\omega_q)^{1/(p-m)}+O(1), &\nu\leq m<j \,,\\
(\lambda/\omega_q)^{1/(p-m)}(\psi_{\nu}^{-1}-\psi_j^{-1})+O(1),  &\nu>m,j>m\,.
\end{cases}
\end{equation*}
It  is easy to show  when $\nu>m$,  $a_{\nu}$ is of order $\lambda^{p/(m-p)}$ and  when $1\leq \nu \leq m$, (\ref{uni:a2}) is still valid. Notice that in this case $\omega_q$ is  a constant  that only depends  on $p$. So now we have finished the proof of Proposition~\ref{uni:prop_anu}.

\subsection{Derivation of $\tilde{a}_{k, \nu}$}
\label{sec:uni:a_tilde}
In this subsection, we shall derive the form of $\tilde{a}_{k, \nu}$ satisfying the constraints in~(\ref{uni:constraint2}). Instead of giving a proposition, we derive the form of $\tilde{a}_{k, \nu}$ in the context.

Consider the $k$'s satisfying $k\in (Kx-p-1, Kx+p+1)$. Since $x$ goes to 0 at a rate of $\lambda^{1/(2m)}/K$, $k>(p+m)$. Hence $\left\{\BO{S}_{k}+\BO{R}_k(x)\right\}^T\BO{\Lambda}_k=1$ is automatically satisfied for arbitrary $\tilde{\BO{a}}_k$.  Denote $\BO{P}=\mathbf{ D}^T\mathbf{D}$ and $\BO{P}_k$ the $k$th column of $\BO{P}$. Note that every row of $\mathbf{B}^T\mathbf{B}/M$ sums to 1, hence 
$$
\left\{\BO{S}_{k}+\BO{R}_k(x)\right\}^T(\Bo{\Lambda}_j -\lambda \mathbf{P}_j)= O\{\lambda^{-1/(2m)}\}+O\left(\max_{1\leq\nu\leq q} |\tilde{a}_{k,\nu}|\right), \quad j=1,\dots, q.
$$
In light of the constraints in (\ref{uni:constraint2}),
$$
\left\{\BO{S}_{k}+\BO{R}_k(x)\right\}^T\mathbf{P}_j = O\left\{\lambda^{-1-1/(2m)}\right\}+\lambda^{-1}O\left(\max_{1\leq\nu\leq q} |\tilde{a}_{k,\nu}|\right), \quad j=1,\dots, q.
$$
For simplicity, denote $O\left\{\lambda^{-1-1/(2m)}\right\}+\lambda^{-1}O\left(\max_{1\leq\nu\leq q} |\tilde{a}_{k,\nu}|\right)$ by  $\xi$.  Further simplification shows that the above is equivalent to
\begin{equation}
\label{uni:tildea1}
\sum_{\nu=1}^q (1-\rho_{\nu}^{-1})^{m+j-1} a_{\nu}\rho_{\nu}^{k-1}+ \sum_{\nu=1}^q (1-\rho_{\nu})^{m+j-1}\tilde{a}_{k, \nu}=O(\xi), \quad j=1,\dots, m,
\end{equation}
and if $p>m$,
\begin{equation}
\label{uni:tildea2}
\sum_{\nu=1}^q (1-\rho_{\nu}^{-1})^{2m} \rho_{\nu}^{-(j-m-1)}a_{\nu}\rho_{\nu}^{k-1}+ \sum_{\nu=1}^q (1-\rho_{\nu})^{2m} \rho_{\nu}^{j-m-1}\tilde{a}_{k, \nu}=O(\xi), \quad j=m+1,\dots, q.
\end{equation}
 
 \subsubsection{The case $p\leq m$}\label{anu_pleqm}

Because $k\in (Kx-p-1, Kx+p+1)$,  $k/\{c_x \lambda^{1/(2m)}\}\rightarrow 1$. Hence for $1\leq\nu\leq m$, $\rho_{\nu}^{k-1}\rightarrow \exp(-c_x\psi_{\nu})$. Since $q=m$,  all $\rho_{\nu}$'s take the forms in (\ref{uni:rho}).  As $\lambda\rightarrow\infty$, $\rho_{\nu}\rightarrow 1$,  $(1-\rho_{\nu})^j \rightarrow \psi_{\nu}^j \lambda^{-j/(2m)}$,  $(1-\rho_{\nu}^{-1})^j \rightarrow (-1)^j \psi_{\nu}^j \lambda^{-j/(2m)}$ and $a_{\nu} \rightarrow \frac{1}{2m}\psi_{\nu} \lambda^{-1/(2m)}$. It is easy to show the leading term of $\sum_{\nu=1}^m (1-\rho_{\nu}^{-1})^{m+j-1} a_{\nu}\rho_{\nu}^{k-1}$ is $ (2m)^{-1}\lambda^{-(m+j)/(2m)}\sum_{\nu=1}^m (-1)^{m+j-1} \psi_{\nu}^{m+j} \exp(-c_x\psi_{\nu})$ and the leading term of $\sum_{\nu=1}^m (1-\rho_{\nu})^{m+j-1}\tilde{a}_{k, \nu}$ is $\lambda^{-(m+j-1)/(2m)} \sum_{\nu=1}^m \psi_{\nu}^{m+j-1} \tilde{a}_{k,\nu}$. Therefore,  we derive that
\begin{equation}
\label{uni:tildea31}
\tilde{a}_{k, \nu}=\frac{\tilde{b}_{k, \nu}}{2m} \lambda^{-1/(2m)}+O(\lambda^{-1/m}), \quad1\leq\nu\leq m,
\end{equation}
for some constant $\tilde{b}_{k, \nu}$. Because of (\ref{uni:tildea31}), $\xi = O\{\lambda^{-1-1/(2m)}\}$. Matching the coefficients of $\lambda^{-(m+j)/(2m)}$ for the $j$th term in  (\ref{uni:tildea1}) gives
\begin{equation}
\label{uni:tildea3}
\sum_{\nu=1}^m (-1)^{m+j-1} \psi_{\nu}^{m+j} \exp(-c_x\psi_{\nu})+ \sum_{\nu=1}^m \psi_{\nu}^{m+j-1} \tilde{b}_{k,\nu} = 0
\end{equation}
To simplify notation, we define $\BO{\Psi}_{m,1}$ is an $m\times m$ matrix with its $(i,j)$th element $\psi_j^{m+i-1}$, $\BO{\Psi}_{m,2}$ is an $m\times m$ matrix with its  $(i,j)$th element $(-1)^{m+j}\psi_j^{m+i}$ and  $\BO{r}(x) = (e^{-\psi_1 x},\dots,e^{-\psi_m x})^T$. By~(\ref{uni:tildea3}), 
\begin{align}
\label{uni:tildeb}
(\tilde{b}_{k, 1},\dots, \tilde{b}_{k, m})^T=\BO{\Psi}_{m,1}^{-1}\BO{\Psi}_{m,2}\BO{r}(c_x).
\end{align}

 \subsubsection{The case $p> m$}
Note that if $\nu>m$, $\rho_{\nu} = O\{\lambda^{-1/(p-m)}\}$ and $a_{\nu} = O\{\lambda^{-p/(p-m)}\}$. Equality~(\ref{uni:tildea2}) for $j=m+1$ reduces to
$$
(-1)^{m+1}\lambda^{-1-1/(2m)}\sum_{\nu=1}^m\psi_{\nu}\exp(-c_x\psi_{\nu}) +(-1)^{m+1} \lambda^{-1} \sum_{\nu=1}^m \tilde{a}_{k,\nu} +\sum_{\nu=m+1}^q \tilde{a}_{k,\nu} = O(\xi),
$$
i.e.,
\begin{equation}
\label{uni:tildea4}
\sum_{\nu=m+1}^q \tilde{a}_{k,\nu} = \lambda^{-1}(-1)^{m+1} \sum_{\nu=1}^m \tilde{a}_{k,\nu} +O(\xi) = O(\xi).
\end{equation}
Because of~(\ref{uni:tildea4}), the analysis in the previous subsection is also valid and~(\ref{uni:tildeb}) still holds. Furthermore, we can derive from~(\ref{uni:tildea2}) that
\begin{equation}
\label{uni:tildea5}
\sum_{\nu=m+1}^q \tilde a_{k,\nu}\rho_{\nu}^j = O\left\{\lambda^{-1-1/(2m)}\right\},\quad j = 0,\dots, q-m-1.
\end{equation}
It follows from~(\ref{uni:tildea5}) that
\begin{equation}
\label{uni:tildeab2}
\sum_{\nu=m+1}^q \tilde a_{k,\nu}\rho_{\nu}^j = O\left\{\lambda^{-1-1/(2m)}\right\},\,\,\text{for any non-negative integer}\,\, j.
\end{equation}

\section{Derivation of Asymptotics}\label{sec:uni:psplines_asymptotics}
In this section, we shall prove the main results in Section~\ref{sec:uni:main}. Specifically, we shall derive  the asymptotic distribution of  P-splines when $x\in (0,1)$ and when $x$ goes to 0 at certain rate. Define $\bar{x}_k = (k-1/2)/K$. 

\subsection{The Case $x\in (0,1)$}\label{subsection:uni:interior}
To prove Proposition~\ref{uni:prop_est}, we need  Proposition~\ref{uni:prop_kernel} below.
\begin{prop}\label{uni:prop_kernel}
Let $h_n = \lambda^{1/(2m)}/K$.  Let $\psi_0 = \min\{\text{Re}(\psi_1),\dots, \text{Re}(\psi_m)\}$, where $\text{Re}(\cdot)$ gives the real part of a complex number. Assume $h_n = o(1)$ and $(Kh_n)^{-1} = o(1)$.
For $x\in (0,1)$, 
\begin{equation}
\label{uni:kernel_H}
\begin{split}
&nh_n \sum_{k,r}B_k(x)B_r(x_i)S_{k,r}/M \\
=&H_m\left(\frac{|x-x_i|}{h_n}\right)+\delta_{\{p>m\}}\left[O\left(\lambda^{-2+\frac{1}{2m}}\right)+ \delta_{\{|x-x_i|< (3p+2-m)/K\}}O\left(\lambda^{-\frac{p}{p-m}+\frac{1}{2m}}\right)\right]\\
&+ \exp\left(-\psi_0 \frac{|x-x_i|}{h_n}\right)\left[ O\left(\lambda^{-1/m}\right)
+ \delta_{\{m=1\}}\delta_{\left \{|x-x_i|\leq (p+1)\lambda^{-1/(2m)}\right \}}O\left\{\lambda^{-1/(2m)}\right\}\right].
\end{split}
\end{equation}
Here $\delta_{\{p>m\}}=1$ if $p>m$ and 0 otherwise; the other $\delta$ terms are similarly defined.
\end{prop}
{\it Proof of Proposition~\ref{uni:prop_kernel}:} By the definition of $\mathbf{S}_{k}$ in~(\ref{uni:S_k}),  
\[
\sum_{k,r}B_k(x)B_r(x_i)S_{k,r}/M =\sum_{\nu=1}^q\left\{\sum_{k,r} B_k(x) B_r(x_i) a_{\nu} \rho_{\nu}^{|k-r|}/M\right\}.
\]
If $p>m$ and  $\nu>m$, $\rho_{\nu} = O\{\lambda^{-1/(p-m)}\}$ by Proposition~\ref{uni:prop_rho_1} and $a_{\nu}$ is of order $\lambda^{-p/(p-m)}$ by Proposition~\ref{uni:prop_anu}.  Note that if $|x-x_i|\geq (3p+2-m)/K$,  a necessary condition for a nonzero $B_k(x)B_r(x_i)$ is that $|k-r|\geq p-m$, hence, for $\nu>m$,
\begin{equation}
 \label{uni:kernel_H1}
 \begin{split}
&\sum_{k,r}B_k(x)B_r(x_i)a_{\nu}\rho_{\nu}^{|k-r|}/M \\
=&\delta_{\{|x-x_i|< (3p+2-m)/K\}}O\left\{\lambda^{-p/(p-m)}Kn^{-1}\right\} +O(\lambda^{-2}Kn^{-1}).
\end{split}
\end{equation}
In the above derivation, Lemma~\ref{uni:lem1} was used. Fix $1\leq \nu\leq m$. Define 
\begin{equation*}
b_{\nu} = -\lambda^{1/(2m)} \log(\rho_{\nu}), \quad 1\leq \nu \leq m.
\end{equation*} 
Then by (\ref{uni:rho}),
\begin{equation*}
\label{uni:bnu}
b_{\nu} = \psi_{\nu} +O\left(\lambda^{-1/m}\right), \quad 1\leq \nu\leq m.
\end{equation*}
It follows that
\begin{equation*}
\rho_{\nu}^{|k-r|} = \exp\left(-b_{\nu} \frac{|\bar{x}_k-\bar{x}_r|}{h_n}\right)= \exp\left(-\psi_{\nu}  \frac{|\bar{x}_k-\bar{x}_r|}{h_n}\right)\left \{ 1+ \frac{|\bar{x}_k-\bar{x}_r|}{h_n}O\left(\lambda^{-1/m}\right)\right\}.
\end{equation*}
By the expression of $a_{\nu}$ in~(\ref{uni:a2}), 
\begin{equation*}
a_{\nu}\rho_{\nu}^{|k-r|} = \frac{\psi_{\nu}}{2mKh_n}\exp\left(-\psi_{\nu}  \frac{|\bar{x}_k-\bar{x}_r|}{h_n}\right)\left \{ 1+ \left(1+\frac{|\bar{x}_k-\bar{x}_r|}{h_n}\right)O\left(\lambda^{-1/m}\right)\right\}.
\end{equation*}
In light of Lemma~\ref{uni:lem3}, 
\begin{align}
\label{uni:kernel_H2}
&2mnh_n\left\{\sum_{k,r} B_k(x) B_r(x_i) a_{\nu} \rho_{\nu}^{|k-r|}/M\right\}\nonumber\\
=&\sum_{k,r} B_k(x) B_r(x_i) \psi_{\nu}\exp\left(-\psi_{\nu}  \frac{|\bar{x}_k-\bar{x}_r|}{h_n}\right)\left \{ 1+ \left(1+\frac{|\bar{x}_k-\bar{x}_r|}{h_n}\right)O\left(\lambda^{-1/m}\right)\right\}\nonumber\\
=&\psi_{\nu}\exp\left(-\psi_{\nu}  \frac{|x-x_i|}{h_n}\right)\left \{ 1-\frac{\psi_{\nu}}{Kh_n}\tilde{g}(x,x_i)+ O\left(\lambda^{-1/m}\right)\right\}.
\end{align}
Summing~(\ref{uni:kernel_H2}) for $\nu=1,\dots, m$ gives
\begin{align}
&nh_n\left\{\sum_{\nu=1}^m\sum_{k,r} B_k(x) B_r(x_i) a_{\nu} \rho_{\nu}^{|k-r|}/M\right\}\nonumber\\
=&\frac{1}{2m}\sum_{\nu=1}^m \psi_{\nu}\exp\left(-\psi_{\nu}  \frac{|\bar{x}-x_i|}{h_n}\right)\left \{ 1-\frac{\psi_{\nu}}{Kh_n}\tilde{g}(x,x_i)+ O\left(\lambda^{-1/m}\right)\right\}\nonumber\\
=&H_m\left(\frac{|x-x_i|}{h_n}\right) + \exp\left(-\psi_0 \frac{|x-x_i|}{h_n}\right) O\left(\lambda^{-1/m}\right)-\frac{1}{Kh_n} \tilde{g}(x,x_i)Q\left(\frac{|x-x_i|}{h_n}\right),
\end{align}
where 
\[
Q(x) = \frac{1}{2m} \sum_{\nu=1}^m \psi_{\nu}^2 \exp\left(-\psi_{\nu}|x|\right).
\]
It is easy to show that $|Q(x)|\leq \exp(-\psi_0 |x|)$.  Lemma~\ref{uni:lem4} states that $\tilde{g}(x,x_i)=0$ if $|x-x_i|\geq (p+1)/K$. Lemma~\ref{uni:lem8} states when $m>1$, $\sum_{1\leq \nu \leq m} \psi_{\nu}^2 =0 $. Thus  if $x$ is close to 0 and $m>1$,
$\sum_{1\leq\nu\leq m} \psi_{\nu}^2 \exp(-\psi_{\nu}|x|) $ is of  the same order as $x$. Hence, 
\begin{equation}
\label{uni:kernel_H3}
\begin{split}
&\tilde{g}(x,x_i)Q\left(\frac{|x-x_i|}{h_n}\right)\\
=&\delta_{\left\{|x-x_i|\leq (p+1)/(Kh_n)\right\}}\exp\left(-\psi_0  \frac{|x-x_i|}{h_n}\right)\left[O\left\{(Kh_n)^{-2}\right\} + \delta_{\{m=1\}} O\left\{(Kh_n)^{-1}\right\}\right].
\end{split}
\end{equation} 
Equalities~(\ref{uni:kernel_H1})--(\ref{uni:kernel_H3}) together prove Proposition~\ref{uni:prop_kernel}. 

{\it Proof of Proposition~\ref{uni:prop_est}:} By~(\ref{uni:est5}) and Proposition~\ref{uni:prop_kernel},
$$
\hat{\mu}(x) = \frac{1}{nh_n} \sum_{i=1}^n y_i\left\{H_m \left(\frac{|x-x_i|}{h_n}\right)+r_i(x)\right\} = \mu^{\ast}(x) + \frac{1}{nh_n}\sum_{i=1}^n r_i(x) y_i,
$$
where 
\begin{equation}
\label{uni:r_i}
\begin{split}
r_i(x) =& \exp\left(-\psi_0 \frac{|x-x_i|}{h_n}\right)\left[ O\left(\lambda^{-\frac{1}{m}}\right)
+ \delta_{\{m=1\}}\delta_{\left \{|x-x_i|\leq (p+1)\lambda^{-1/(2m)}\right\}}O\left(\lambda^{-\frac{1}{2m}}\right)\right]\\
&+\delta_{(p>m)}\left[O\left(\lambda^{-2+\frac{1}{2m}}\right)+\delta_{\left\{|x-x_i|< (3p+2-m)/K\right\}}O\left\{\lambda^{-\frac{p}{p-m}+\frac{1}{2m}}\right\}\right]\\
& +O\left[nh_n\exp\{-C \lambda^{-\frac{1}{2m}}K\min(x,1-x)\}\right].
\end{split}
\end{equation}
First we have
\begin{equation}
\label{uni:mu_exp_1}
\left|\textrm{E}\left\{\hat{\mu}(x)-\mu^{\ast}(x)\right\} \right|\leq (nh_n)^{-1}\sum_i \left|\mu(x_i)r_i(x)\right|.
\end{equation}
We study the right hand side of~(\ref{uni:mu_exp_1}). For $r_i(x)$ defined in~(\ref{uni:r_i}),  the two terms $O\{\lambda^{-2+1/(2m)}\}$ and $O\left[nh_n\exp\{-C \lambda^{-1/(2m)}K\min(x,1-x)\}\right]$ are of order $o(\lambda^{-1/m})$. Also
\begin{align*}
&(nh_n)^{-1}\sum_i |\mu(x_i)| \exp \left(-\psi_0 \frac{|x-x_i|}{h_n}\right) = O(1),\\
&(nh_n)^{-1}\sum_i |\mu(x_i)| \exp \left(-\psi_0 \frac{|x-x_i|}{h_n}\right) \delta_{\{|x-x_i|\leq (p+1)\lambda^{-1/(2m)}\}} = O\left(\lambda^{-\frac{1}{2m}}\right),\\
&(nh_n)^{-1}\sum_i |\mu(x_i)| \delta_{\{|x-x_i|\leq (3p+2-m)/K\}}= O\{(Kh_n)^{-1}\}.
\end{align*}
It follows that
$\sum_i \left|\mu(x_i) r_i(x)\right | = O(\lambda^{-1/m}).
$
Next we derive that
\begin{equation}
\label{uni:mu_var_1}
\textrm{var}\left\{\hat{\mu}(x)-\mu^{\ast}(x)\right\} = (nh_n)^{-2}\sum_i r_i^2(x)\sigma^2(x_i).
\end{equation}
With similar derivation as before, we can establish that $(nh_n)^{-1}\sum_i r_i^2(x)\sigma^2(x_i) = o(1)$. Therefore the proposition is proved.

\begin{example} Consider the case $m=2$. Denote the imaginary number by $\imath$. Then $\psi_1 =\frac{1+\imath}{\sqrt{2}}$ and $\psi_2=\frac{1-\imath}{\sqrt{2}}$. Hence the equivalent kernel for $x\in (0,1)$ is
\begin{align*}
\frac{1}{2\sqrt{2}} e^{-\frac{|x-\tilde{x}|}{\sqrt{2}}}\left\{\cos\left(\frac{|x-\tilde{x}|}{\sqrt{2}}\right)+\sin\left( \frac{|x-\tilde{x}|}{\sqrt{2}}\right)\right\}.
\end{align*}
\end{example}
\begin{example} Consider the case $m=3$. Then $\psi_1 = 1, \psi_2 = \frac{1+\sqrt{3}\imath}{2}, \psi_3 = \frac{1-\sqrt{3}\imath}{2}$. Hence the equivalent kernel for $x\in (0,1)$ is
\begin{equation*}
\frac{1}{6}e^{-|x-\tilde{x}|}+\frac{1}{6}e^{-\frac{|x-\tilde{x}|}{2}}\left\{\cos\left(\frac{\sqrt{3}|x-\tilde{x}|}{2}\right)+\sqrt{3}\sin\left(\frac{\sqrt{3}|x-\tilde{x}|}{2}\right)\right\}.
\end{equation*}
\end{example}

{\it Proof of Theorem~\ref{uni:thm1}:} Proposition~\ref{uni:prop_est} shows that the P-spline estimator is asymptotically equivalent to a kernel regression estimator with the kernel function $H_m(x)$. Hence a standard analysis of the kernel regression estimator as in Wand and Jones (1995)  with the kernel function $H_m(x)$  should give us the desired result. The detailed derivation is as follows. First,
\begin{equation*}
\textrm{E}\{\mu^{\ast}(x)\} =  \mu(x) +(-1)^{m+1} h_n^{2m}\mu^{(2m)}(x)+ o(h_n^{2m})
\end{equation*}
and 
\begin{equation*}
\begin{split}
\textrm{var}\left\{\mu^{\ast}(x)\right\} &= \sum_i \sigma^2(x_i) \frac{1}{(nh_n)^2} H_m^2 \left(\frac{|x-x_i|}{h_n}\right)\\
&=\frac{1}{nh_n} \sigma^2(x) \int_{-\infty}^{\infty} H_m^2(s) ds + o\{(nh_n)^{-1}\}.
\end{split}
\end{equation*}
By Proposition~\ref{uni:prop_est}, we obtain
\begin{equation*}
\begin{split}
\textrm{E}\{\hat\mu(x)\} &=  \mu(x) +(-1)^{m+1} h_n^{2m}\mu^{(2m)}(x)+ o(h_n^{2m}) + O\{(nh_n)^{-1}\},\\
\textrm{var}\left\{\hat{\mu}(x)\right\} &= \frac{1}{nh_n} \sigma^2(x) \int_{-\infty}^{\infty} H_m^2(s) ds + o\{(nh_n)^{-1}\},
\end{split}
\end{equation*}
and the proof is straightforward by verifying that  $h_n^{4m}$ and $(nh_n)^{-1}$ are of the same order and $\lambda^{-1/m}= o(h_n^{2m})$.

\subsection{The Boundary Case}\label{subsection:uni:boundary}
By (\ref{uni:est6}) and the derivation in Section~\ref{sec:uni:a_tilde},  we have
\begin{align}
\hat{\mu}(x)&= \frac{1}{M}\sum_{i=1}^n y_i \left [\sum_{k,r} B_k(x)B_r(x_i)\left\{S_{k,r}+R_{k, r}(x)\right\}+b_{i,0}(x)\right]\nonumber\\
&=\frac{1}{M}\sum_{i=1}^n y_i\left\{\sum_{k,r} B_k(x)B_r(x_i)S_{k,r}+b_{i,0}(x)\right\}\label{uni:est61}\\
&\quad+\frac{1}{M}\sum_{i=1}^n y_i\left\{\sum_{k, r}B_k(x)B_r(x_i)R_{k, r}(x)\right\}.\label{uni:est612}
\end{align}
Note that $b_{i,0}(x) = O[\exp\{-C_0 \lambda^{-1/(2m)}K\}]$. The sum in (\ref{uni:est61}) can be similarly analyzed as in Section~\ref{subsection:uni:interior} and we have
\begin{equation*}
\begin{split}
&\frac{1}{M}\sum_{i=1}^n y_i\left\{\sum_{k,r} B_k(x)B_r(x_i)S_{k,r}+b_{i,0}(x)\right\}\\
 =& \frac{1}{nh_n} \sum_{i=1}^ny_i\left[ H_m\left(\frac{|x-x_i|}{h_n}\right)+\exp\left (-\psi_0 \frac{|x-x_i|}{h_n} \right)O\left\{(Kh_n)^{-1}\right\}\right]
\end{split}
\end{equation*}
Now we focus on the second sum (denoted by $\hat{\mu}_b(x)$) in (\ref{uni:est612}). Note that $R_{k, r}(x) = \sum_{\nu=1}^q \tilde{a}_{k, \nu}\rho_{\nu}^{r-1}$. Note also  if $\nu>m$, $\rho_{\nu}= O\{\lambda^{-1/(p-m)}\}$ and~(\ref{uni:tildeab2}) holds. Hence, 
\begin{align*}
\hat{\mu}_b(x) = \frac{1}{2mnh_n}\sum_{i=1}^n y_i\left[ \sum_{\nu=1}^m\sum_{r=1}^c\sum_{k=1}^c B_r(x_i)B_k(x) \tilde{b}_{k, \nu}\rho_{\nu}^{r-1}+O\{(Kh_n)^{-2}\}\right].
\end{align*}
By a similar analysis as in Section~\ref{subsection:uni:interior}, we obtain, aided by Lemma~\ref{uni:lem5}, that
\begin{equation*}
\begin{split}
\hat{\mu}_b(x)& =\frac{1}{2mnh_n}\sum_{i=1}^n y_i\left[ \BO{r}^T(\frac{x_i}{h_n})\BO{\Psi}_{m,1}^{-1}\BO{\Psi}_{m,2}\BO{r}(c_x)+O\{(Kh_n)^{-2}\}\right]\\
&=\frac{1}{2mnh_n}\sum_{i=1}^n y_i\left[ \BO{r}^T(\frac{x_i}{h_n})\BO{\Psi}_{m,1}^{-1}\BO{\Psi}_{m,2}\BO{r}(\frac{x}{h_n})+O\{(Kh_n)^{-2}\}\right].
\end{split}
\end{equation*}
Note that $\BO{\Psi}_{m,1}$, $\BO{\Psi}_{m,2}$ and $\mathbf{r}(x)$ are defined in Section~\ref{sec:uni:a_tilde}. 
In the above derivation, we used the assumption that   $ x/h_n$ converges to $c_x$; we also used (\ref{uni:tildeb}). We define  the equivalent kernel for $\hat{\mu}_b(x)$ as
\begin{align}
\label{uni:kernelb2}
H_{b,m}(x,\tilde{x})=\frac{1}{2m} \BO{r}(\tilde{x})^T\BO{\Psi}_{m,1}^{-1}\BO{\Psi}_{m,2}\BO{r}(x).
\end{align} 
Now we have
\begin{equation}
\label{uni:est62}
\hat{\mu}(x) = \frac{1}{nh_n}\sum_{i=1}^n y_i \left[H_m\left(\frac{|x-x_i|}{h_n}\right)+ H_{b,m}\left(\frac{x}{h_n},\frac{x_i}{h_n}\right)+\exp\left (-\psi_0 \frac{x_i}{h_n} \right)O\left(\frac{1}{Kh_n}\right)\right].
\end{equation}
The above equality shows that when $x$ is near 0, a P-spline estimator is a kernel regression estimator with the equivalent kernel
\begin{equation}
\label{uni:kernelb3_1}
H_m(|x-\tilde{x}|)+H_{b,m}(x,\tilde{x}).
\end{equation}
Next we provide two specific examples of~(\ref{uni:kernelb3_1}).
\begin{example}  Consider the case $m=2$.  It can be shown that
\[
\Psi_{m,1} = \left(\begin{array}{cc}
\imath&-\imath\\
\frac{-1+\imath}{\sqrt{2}}&\frac{-1-\imath}{\sqrt{2}}
\end{array}\right),\quad
\Psi_{m,2}=\left(\begin{array}{cc}
-\frac{-1+\imath}{\sqrt{2}}&-\frac{-1-\imath}{\sqrt{2}}\\
-1&-1
\end{array}\right),
\]
and 
\[
\BO{r}(x) =e^{-\frac{x}{\sqrt{2}}}\left(\begin{array}{c} cos\left(\frac{x}{\sqrt{2}}\right)-\imath \sin\left(\frac{x}{\sqrt{2}}\right)\\
  \cos\left(\frac{x}{\sqrt{2}}\right)+\imath \sin\left(\frac{x}{\sqrt{2}}\right)
  \end{array}\right).
\]
Hence,
\[
H_{b,2}(x,\tilde{x}) = \frac{\sqrt{2}}{4}e^{-\frac{|x+\tilde{x}|}{\sqrt{2}}}\left\{\cos\left(\frac{|x-\tilde{x}|}{\sqrt{2}}\right)+2\cos\left(\frac{x}{\sqrt{2}}\right)\cos\left(\frac{\tilde{x}}{\sqrt{2}}\right)-\sin\left(\frac{x+\tilde{x}}{\sqrt{2}}\right)\right\}.
\]
It follows that the equivalent kernel for $x$ near 0 is
\begin{align*}
\label{kernelmb=2}
&\frac{\sqrt{2}}{4} e^{-\frac{|x-\tilde{x}|}{\sqrt{2}}}\left\{\cos\left(\frac{|x-\tilde{x}|}{\sqrt{2}}\right)+\sin\left( \frac{|x-\tilde{x}|}{\sqrt{2}}\right)\right\}\nonumber\\
+&\frac{\sqrt{2}}{4}e^{-\frac{|x+\tilde{x}|}{\sqrt{2}}}\left\{\cos\left(\frac{|x-\tilde{x}|}{\sqrt{2}}\right)+2\cos\left(\frac{x}{\sqrt{2}}\right)\cos\left(\frac{\tilde{x}}{\sqrt{2}}\right)-\sin\left(\frac{x+\tilde{x}}{\sqrt{2}}\right)\right\}.
\end{align*}
When $x=0$, the equivalent kernel becomes
\[
\sqrt{2}e^{-\tilde{x}/\sqrt{2}}\cos\left(\tilde{x}/\sqrt{2}\right),
\]
which coincides with the equivalent kernel for the smoothing splines (Silverman, 1984).
\end{example}
\begin{example} Consider the case $m=3$. It can be shown that
\[
\Psi_{m,1} = \left(\begin{array}{ccc}
1&-1&-1\\
1&\frac{-1-\sqrt{3}\imath}{2}&\frac{-1+\sqrt{3}\imath}{2}\\
1&\frac{1-\sqrt{3}\imath}{2}&\frac{1+\sqrt{3}\imath}{2}
\end{array}\right),\quad
\Psi_{m,2}=\left(\begin{array}{ccc}
1&\frac{-1-\sqrt{3}\imath}{2}&\frac{-1+\sqrt{3}\imath}{2}\\
1&\frac{1-\sqrt{3}\imath}{2}&\frac{1+\sqrt{3}\imath}{2}\\
1&1&1
\end{array}\right),
\]
and 
\[
\BO{r}(x) =\left(\begin{array}{c} e^{-x}\\e^{-\frac{x}{2}}\left\{cos\left(\frac{\sqrt{3}x}{2}\right)-\imath \sin\left(\frac{\sqrt{3}x}{2}\right)\right\}\\
  e^{-\frac{x}{2}}\left\{\cos\left(\frac{\sqrt{3}x}{2}\right)+\imath \sin\left(\frac{\sqrt{3}x}{2}\right)\right\}
  \end{array}\right).
\]
It follows that the equivalent kernel for $x$ near 0 is
\begin{align*}
\label{kernelmb=3}
&\frac{1}{6}e^{-|x-\tilde{x}|}+\frac{1}{6}e^{-\frac{|x-\tilde{x}|}{2}}\left\{\cos\left(\frac{\sqrt{3}|x-\tilde{x}|}{2}\right)+\sqrt{3}\sin\left(\frac{\sqrt{3}|x-\tilde{x}|}{2}\right)\right\}\nonumber\\
+&\frac{3}{6}e^{-|x+\tilde{x}|} +\frac{2}{6}e^{-|x+\frac{\tilde{x}}{2}|}\left\{\cos\left(\frac{\sqrt{3}\tilde{x}}{2}\right)-\sqrt{3}\sin\left(\frac{\sqrt{3}\tilde{x}}{2}\right)\right\}\\
+&\frac{2}{6}e^{-|\tilde{x}+\frac{x}{2}|}\left\{\cos\left(\frac{\sqrt{3}x}{2}\right)-\sqrt{3}\sin\left(\frac{\sqrt{3}x}{2}\right)\right\}\\
+&\frac{1}{6}e^{-\frac{|x+\tilde{x}|}{2}}\left\{3\cos\left(\frac{\sqrt{3}(\tilde{x}-x)}{2}\right)-\sqrt{3}\sin\left(\frac{\sqrt{3}(\tilde{x}-x)}{2}\right)+2\sin\left(\frac{\sqrt{3}x}{2}\right)\sin\left(\frac{\sqrt{3}\tilde{x}}{2}\right)\right\}.
\end{align*}
When $x=0$, the equivalent kernel becomes
\[
e^{-\tilde{x}}+e^{-\tilde{x}/2}\left\{ \cos\left(\frac{\sqrt{3}\tilde{x}}{2}\right)-\frac{\sqrt{3}}{3}\sin\left(\frac{\sqrt{3}\tilde{x}}{2}\right)\right\}.
\]
\end{example}

{\it Proof of Theorem~\ref{uni:thm2}:} Similar to the proof of Theorem~\ref{uni:thm1}, we can derive that
\begin{equation*}
\begin{split}
&\textrm{E}\{\hat{\mu}(x)\}\\
 =& \frac{1}{nh_n}\sum_{i=1}^n \mu(x_i) \left[H_m\left(\frac{|x-x_i|}{h_n}\right)+ H_{b,m}\left(\frac{x}{h_n},\frac{x_i}{h_n}\right)+\exp\left (-\psi_0 \frac{x_i}{h_n} \right)O\left(\frac{1}{Kh_n}\right)\right]\\
=&\frac{1}{h_n}\int_0^1 \mu(u) \left\{H_m\left(\frac{|x-u|}{h_n}\right)+ H_{b,m}\left(\frac{x}{h_n},\frac{u}{h_n}\right)\right\}du+O\left(\frac{1}{Kh_n}\right)\\
=&\int^{c_x}_{-\infty} \mu(x-hv)\left\{H_m(v)+H_{b,m}(c_x,c_x-v)\right\}dv +O\left\{(Kh_n)^{-1}\right\},
\end{split}
\end{equation*}
and
\begin{equation}
\label{uni:var_mu_b}
\begin{split}
&\textrm{var}\{\hat{\mu}(x)\}\\
 =& \frac{1}{(nh_n)^2}\sum_{i=1}^n \sigma^2(x_i) \left[H_m\left(\frac{|x-x_i|}{h_n}\right)+ H_{b,m}\left(\frac{x}{h_n},\frac{x_i}{h_n}\right)+\exp\left (-\psi_0 \frac{x_i}{h_n} \right)O\left(\frac{1}{Kh_n}\right)\right]^2\\
=&\frac{1+o(1)}{nh_n}\frac{1}{h_n}\int_0^1 \sigma^2(u) \left\{H_m\left(\frac{|x-u|}{h_n}\right)+ H_{b,m}\left(\frac{x}{h_n},\frac{u}{h_n}\right)\right\}^2du\\
=&\frac{1+o(1)}{nh_n} \sigma^2(x)\int_{-\infty}^{c_x} \left\{H_m(v)+H_{b,m}(c_x, c_x-v)\right\}^2 dv.
\end{split}
\end{equation}
By Proposition~\ref{uni:prop_b} below, we have
\begin{equation}
\label{uni:mu_exp_b}
\begin{split}
\textrm{E}\{\hat{\mu}(x)\} &=\mu(x) + (-1)^{m+1}h_n^{m}\mu^{(m)}(x)\int^{c_x}_{-\infty} v^m\left\{H_m(v)+H_{b,m}(c_x, c_x-v)\right\}dv \\
&\quad+ o\left(h_n^{m+1}\right)+O\left\{(Kh_n)^{-1}\right\}.
\end{split}
\end{equation}
Combining~(\ref{uni:var_mu_b}) with~(\ref{uni:mu_exp_b}), Theorem~\ref{uni:thm2} is proved.

\begin{prop}\label{uni:prop_b}
For any fixed constant $t\geq 0$, 
$$\int^{t}_{-\infty} x^{\ell} \left\{H_m(x)+H_{b,m}(t,t-x)\right\}dx  = 0, \quad \ell = 1,\dots, m-1,$$
and 
$$\int^{t}_{-\infty} x^m\left\{H_m(x)+H_{b,m}(t,t-x)\right\}dx \neq 0.$$
\end{prop}
{\it Proof of Proposition~\ref{uni:prop_b}:}
By Lemma~\ref{uni:lem6}, we can show that
\begin{equation*}
\begin{split}
\int^{t}_{-\infty} x^{\ell} H_m(x)\mathrm{d}x &=-\frac{\ell!}{2m}\sum_{k=1}^{\ell+1}\sum_{\nu=1}^m \frac{t^{\ell-k+1}}{(\ell-k+1)|}\bar{\psi}_{\nu}^{k-1}e^{-\psi_{\nu} t}\\
&=-\frac{\ell!}{2m}\left\{\sum_{k=1}^{\ell+1}\frac{t^{\ell-k+1}}{(\ell-k+1)!}\bar{\psi}_1^{k-1},\dots, \sum_{k=1}^{\ell+1}\frac{t^{\ell-k+1}}{(\ell-k+1)!}\bar{\psi}_m^{k-1} \right\}\mathbf{r}(t),
\end{split}
\end{equation*}
and
\begin{equation*}
\int_{-\infty}^t x^{\ell} \mathbf{r}(t-x)^T\mathrm{d} x= -\ell! \left\{\sum_{k=1}^{\ell+1} \frac{t^{\ell-k+1}}{(\ell-k+1)!} (-1)^k \bar{\psi}_1^k,\dots, \sum_{k=1}^{\ell+1} \frac{t^{\ell-k+1}}{(\ell-k+1)!} (-1)^k \bar{\psi}_m^k\right\}.
\end{equation*}
Because  $H_{b,m}(t, t-x) = (2m)^{-1}\mathbf{r}(t-x)^T \BO{\Psi}_{m,1}^{-1}\BO{\Psi}_{m,2}\mathbf{r}(t)$, it suffices to prove that
\begin{equation}
\label{uni:eqn_b}
\left(\bar{\psi}_1^{k-1},\dots, \bar{\psi}_m^{k-1} \right)+(-1)^k \left(\bar{\psi}_1^k, \dots, \bar{\psi}_m^k\right)\BO{\Psi}_{m,1}^{-1}\BO{\Psi}_{m,2} = \mathbf{0}^T,\quad k=1,\dots, m.
\end{equation}
Let $
\mathbf{w}_k^T = (-1)^{m+1} \left(\bar{\psi}_1^k, \dots, \bar{\psi}_m^k\right)\BO{\Psi}_{m,1}^{-1}\BO{\Psi}_{m,2}$. Then $\mathbf{w}_k$ is the $(m+1-k)$th row of $\mathbf{\Psi}_{m,2}$. Hence,  for $k=1,\dots, m$, 
\[
\mathbf{w}_k^T = (-1)^{2m-k+1}\left( \psi_1^{2m-k+1},\dots, \psi_m^{2m-k+1}\right)=(-1)^{m-k}\left(\bar{\psi}_1^{k-1},\dots, \bar{\psi}_m^{k-1} \right)
\]
which proves~(\ref{uni:eqn_b}). For $\ell = m$, we have
\[
\int^{t}_{-\infty} x^m \left\{H_m(x)+H_{b,m}(t,t-x)\right\}dx =\frac{-m!}{2m} \tilde{\mathbf{w}}_{m+1}^T\mathbf{r}(t),
\]
where $\tilde{\mathbf{w}}_{m+1}^T = \left(\bar{\psi}_1^m,\dots, \bar{\psi}_m^m \right)+ (-1)^{m+1}\left(\bar{\psi}_1^{m+1},\dots, \bar{\psi}_m^{m+1}\right)\BO{\Psi}_{m,1}^{-1}\BO{\Psi}_{m,2}$. Note that $\left(\psi_1^m,\dots, \psi_m^m\right)  = (-1)^{m+1}\left(\bar{\psi}_1^{m+1},\dots, \bar{\psi}_m^{m+1}\right) $ is the first row of $\mathbf{\Psi}_{m,1}$, hence
\begin{equation*}
\begin{split}
\tilde{\mathbf{w}}_{m+1}^T &= \left(\bar{\psi}_1^m,\dots, \bar{\psi}_m^m \right)+(-1)^{m+1}\left(\psi_1^m, \dots, \psi_m^m\right)\\
&=2(-1)^{m+1}\left(\psi_1^m, \dots, \psi_m^m\right)
\end{split}
\end{equation*}
which finishes the proof.

\section{Irregularly Spaced Data}\label{sec:uni:irregular}
Suppose the design points $\underbar{x} = \{x_1,\dots, x_n\}$ are independent and sampled from a distribution $F(x)$ in $[0,1]$. Suppose $F(x)$ is twice continuously differentiable with derivative $f(x)$ and $f(x)$ is positive over $[0,1]$.  For unequally spaced design points,  the asymptotic analysis in Section~\ref{sec:uni:psplines_asymptotics} does not hold here. Instead of pursuing the challenging task of analyzing the P-splines fitted to irregularly spaced data directly, we first bin the data. So we partition $[0,1]$ into $I$ intervals with equal lengths, and let $\tilde{y}_k$ be the mean of all $y_i$ such that $x_i$ is in the $k$th bin. If the $k$th bin has no data point, we let $\tilde{y}_k$ be 0. Here we assume $I\sim c_I n^{\tau_I}$ for some constants $c_I$ and $\tau_I<1$.  Assuming $\tilde{y}_{k}$ is the data point at $\tilde{x}_{k}$, the center of the $k$th bin, we apply  P-splines to the binned data $(\tilde{y}_{k})_{1\leq k\leq I}$ to get
$$
\hat{\boldsymbol{\theta}}^{\ast} = \BO{\Lambda}^{-1}\BO{B}^T\tilde{\BO{y}}/M.
$$
Then the penalized estimate is defined as
\begin{equation}
\label{uni:est_irregular}
\hat{\mu}(x) =\sum_{k=1}^c \hat{\theta}_{k}^{\ast}B_{k}(x).
\end{equation}
Note that the practice of  binning data in penalized splines also appears in Wang and Shen (2010). The asymptotic distribution  of $\hat{\mu}(x)$ in~(\ref{uni:est_irregular}) can be similarly derived as in Section~\ref{sec:uni:psplines_asymptotics}. 
\begin{thm}\label{uni:thm3}
Let $\sigma^2(x) = \textrm{\textnormal{var}}(y|X=x)$. Assume $\tau_I>\max(\tau, 1/2)$ and condition (1)-(4) in Proposition~\ref{uni:prop_est} hold. Furthermore, assume $\sigma^2(x)$ has a continuous second derivative. For $x\in (0,1)$,  with the same notation and assumptions as in Theorem~\ref{uni:thm1}, we have that
\[
n^{2m/(4m+1)}\left\{\hat{\mu}(x)-\mu(x)\right\}\Rightarrow N\left\{\tilde{\mu}(x),V(x)/f(x) \right\}
\]
in distribution as $n\rightarrow\infty$, where $\tilde{\mu}(x)$ is defined  in~(\ref{uni:tilde_mu(x)}) and  $V(x)$ is defined in~(\ref{uni:V(x)}).
 \end{thm}
 \begin{rem}
 The above theorem holds for the fixed design as well and the assumption required for the design points is an analogue to~(\ref{uni:n_k}): $\sup_k \left|n_k/(nI^{-1}) - f(\tilde{x}_{\kappa}) \right| = o(1)$.
 \end{rem}
{\it Proof of Theorem~\ref{uni:thm3}:}  By a similar analysis as in Section~\ref{sec:uni:psplines_asymptotics} to the binned data $\tilde{\BO{y}}$ and with $n$ replaced by $I$, we obtain
$$
\hat{\mu}(x) = \frac{1}{Ih_n}\sum_{k=1}^I \tilde{y}_k\left\{ H_m \left(\frac{|x-\tilde{x}_k|}{h_n}\right)+r_k(x)\right\},
$$
where 
\begin{equation*}
\begin{split}
r_k(x) =& \exp\left(-\psi_0 \frac{|x-\tilde x_k|}{h_n}\right)\left[ O\left(\lambda^{-1/m}\right)
+ \delta_{\{m=1\}}\delta_{\left \{|x-\tilde x_k|\leq (p+1)\lambda^{-1/(2m)}\right \}}O\left\{\lambda^{-1/(2m)}\right\}\right]\\
&+\delta_{(p>m)}\left[O\left(\lambda^{-2+\frac{1}{2m}}\right)+\delta_{\left\{|x-\tilde x_k|< (3p+2-m)/K\right\}}O\left\{\lambda^{-\frac{p}{p-m}+\frac{1}{2m}}\right\}\right]\\
& +O\left[Ih_n\exp\{-C \lambda^{-1/(2m)}K\min(x,1-x)\}\right].
\end{split}
\end{equation*}
Then
\begin{equation}
 \label{uni:mu1}
\textrm{E}\left\{\hat{\mu}(x)|\underbar{x}\right\} = (Ih_n)^{-1}\sum_k \textrm{E}\left\{\tilde{y}_k|\underbar{x}\right\}\left\{H_m\left(\frac{x-\tilde{x}_k}{h_n}\right)+r_k(x)\right \},
\end{equation}
and
\begin{equation}
\label{uni:var1}
\textrm{var}\left\{\hat{\mu}(x)|\underbar{x}\right\}=(Ih_n)^{-2}\sum_{k} \textrm{var}\left\{\tilde{y}_k|\underbar{x}\right\}\left\{H_{m}\left(\frac{x-\tilde x_k}{h_{n}}\right)+r_k(x)\right\}^2.
\end{equation}
For simplicity, we let 
$$
G_k = H_m\left\{h_n^{-1}(x-\tilde{x}_k)\right\} + b_k(x).
$$
Let $n_k$ be the number of data points in the $k$th bin, then 
\[
 \textrm{var}\left\{\tilde{y}_k|\underbar{x}\right\} = n_k^{-2}\sum_{i=1}^n \sigma^2(x_i) \delta_{\{|x_i-\tilde{x}_k|\leq (2I)^{-1}\}}.
\]
So $\textrm{var}\left\{\sqrt{n_k}\tilde{y}_k |\underbar{x} \right\}$ is a Nadaraya-Watson kernel regression estimator of the conditional variance function $\sigma^2(x)$ at $\tilde{x}_k$. Similarly, $n_k/(nI^{-1})$ is a kernel density estimator of $f(x)$ at $\tilde{x}_k$.  By the uniform convergence theory for kernel density estimators and Nadaraya-Watson kernel regression estimators (see, for instance, Hansen (2008)),
\begin{equation}
\label{uni:n_k}
\sup_k \left|n_k/(nI^{-1}) - f(\tilde{x}_{\kappa}) \right| = O_p\left\{\sqrt{I\ln n/n} + I^{-2}\right\} = o_p(1),
\end{equation}
and
\begin{equation*}
\sup_k \left| \textrm{var}\left\{ \sqrt{n_k}\tilde{y}_k |\underbar{x}\right\}-\sigma^2(\tilde{x}_k)\right| = O_p\left\{\sqrt{I\ln n/n} + I^{-2}\right\} = o_p(1).
\end{equation*}
It follows that
\begin{equation}
\label{uni:var2}
\sup_k \left| \frac{n}{I}\textrm{var}\left\{\tilde{y}_k |\underbar{x}\right\}-\frac{\sigma^2(\tilde{x}_{\kappa})}{f(\tilde{x}_{\kappa})}\right| =o_p(1) .
\end{equation}
Then by~(\ref{uni:var1}) and~(\ref{uni:var2}), 
\begin{equation*}
\left|\textrm{var}\left\{\hat{\mu}(x)|\underbar{x}\right\}-\frac{1}{nh_nIh_n}\sum_k \frac{\sigma^2(\tilde{x}_{\kappa})}{f(\tilde{x}_{\kappa})}G_k^2\right|
=\frac{o_p(1)}{nh_nIh_n} \sum_k G_k^2 = o_p\left\{(nh_n)^{-1}\right\},\label{as1}
\end{equation*}
and hence
\begin{equation}
\label{uni:var3}
\text{var}\left\{\hat{\mu}(x)|\underbar{x}\right\} = \frac{1}{nh_n} \frac{V(x)}{f(x)} + o_p\left\{(nh_n)^{-1}\right\}.
\end{equation}
where $V(x)$ is defined in~(\ref{uni:V(x)}). Because
\begin{equation*}
\textrm{E}\left\{\tilde{y}_k|\underbar{x}\right\} =n_k^{-1}\sum_{i=1}^n \mu(x_i) \delta_{\{|x_i-\tilde{x}_{\kappa}|\leq (2I)^{-1}\}},
\end{equation*}
we can derive by~(\ref{uni:n_k})  that
\begin{equation*}
\sup_k \left|\textrm{E}\left\{\tilde{y}_k |\underbar{x}\right\}- \mu(\tilde{x}_{\kappa}) \right| = O_p(I^{-1}).
\end{equation*}
Hence by~(\ref{uni:mu1}),
\begin{equation*}
\left| \textrm{E}\left\{\hat{\mu}(x)|\underbar{x}\right\}- \frac{1}{Ih_n}\sum_k \mu(\tilde{x}_{\kappa})G_k\right| =  O_p(I^{-1}),
\end{equation*}
and hence
\begin{equation}
\label{uni:mu2}
\textrm{E}\left\{\hat{\mu}(x)|\underbar{x}\right\}=\mu(x) + n^{-(2m)/(4m+1)} \tilde{\mu}(x) + o_p\left\{n^{-(2m)/(4m+1)}\right\},
\end{equation}
where $\tilde{\mu}(x)$ is defined in~(\ref{uni:tilde_mu(x)}). 
With~(\ref{uni:var3}) and~(\ref{uni:mu2}), we can derive that
\begin{equation}
\label{uni:as4}
n^{(2m)/(4m+1)}\left[\hat{\mu}(x) -\textrm{E}\left\{\hat{\mu}(x)|\underbar{x}\right\} \right]\Rightarrow N\left\{0, V(x)/f(x)\right\}
\end{equation}
in distribution and
\begin{equation}
\label{uni:as5}
n^{(2m)/(4m+1)}\left[\textrm{E}\left\{\hat{\mu}(x)|\underbar{x}\right\}-\mu(x) \right] = \tilde{\mu}(x) + o_p(1).
\end{equation}
Equalities~(\ref{uni:as4}) and (\ref{uni:as5}) together prove the theorem.


\section{An Example}\label{sec:uni:example}
\begin{figure}[t]
\begin{center}
\includegraphics[scale=0.8,angle=270]{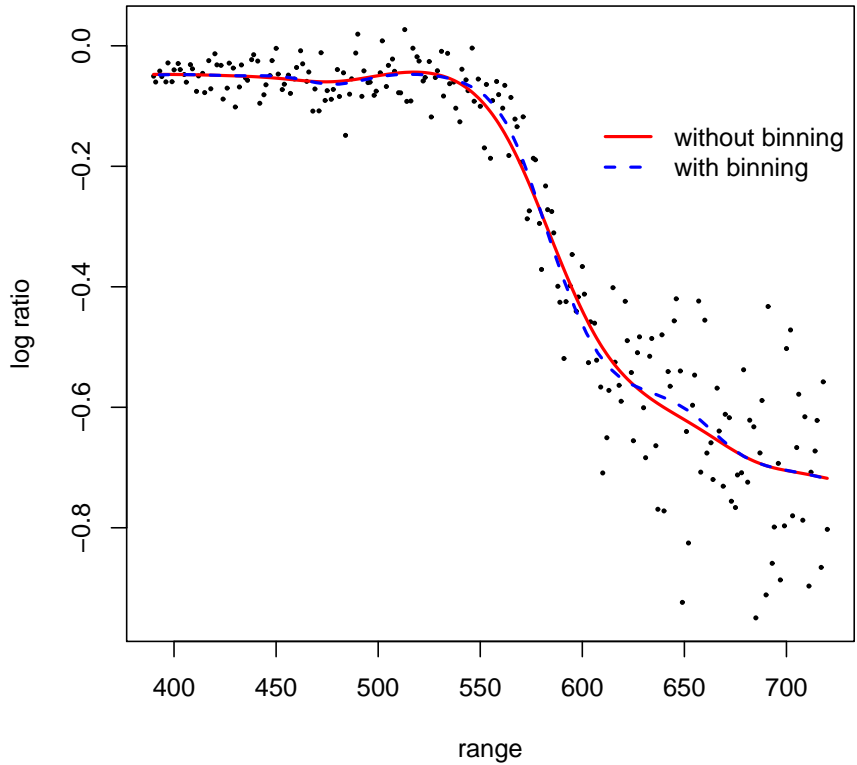}
\caption{\label{uni:lidar} The fitted curves of the response, log ratio,  as a function of the predictor, range. The solid line is the fitted P-splines without binning the data, and the dashed line is the fitted P-splines after binning the data. The solid dots are the observed data. }
\end{center}
\end{figure}

We illustrate the idea of binning data using the LIDAR (light detection and ranging) data. The LIDAR data were analyzed in Holst et al.\ (1996) and Ruppert et al.\ (1997). The LIDAR data have $221$ data points, and  details about the LIDAR data can also be found in Ruppert et al.\ (2003). We fit the response, {\it logratio}, as a function of the predictor, {\it range}. First, we fit the data using cubic P-splines with a penalty of second order, and we use 35 equidistant knots  as suggested in Ruppert et al.\ (2003).  Then,  we fit the binned data using cubic P-splines with a penalty of second order. The number of bins is 60 and we use 15 equidistant knots. The result is given in  Figure~\ref{uni:lidar}. We can see that the two fitted curves are similar, with biggest difference occurring when the predictor, {\it range}, is around 650.

\section{Discussion}\label{sec:uni:conclusion}
We have concentrated on the asymptotics of penalized splines estimation. In contrast to smoothing splines, penalized splines allow us to choose the number of knots, the degree of splines and the penalty independently. Our study provides theoretical guidelines on how to choose them. In our setting, the penalty $\lambda$ plays the role of a smoothing parameter and the optimal order for $\lambda$ is provided. The number of knots $K$ is not important as long as it exceeds a given bound. The choice of the degree of splines does not affect the asymptotic distribution. Our results indicate that the performance of penalized splines estimation is similar to that of smoothing splines estimation (Silverman, 1984) and a class of kernel estimators (Messer and Goldstein, 1993). Furthermore, penalized splines have a slower convergence rate at the boundary than in the interior.

\section{Some Lemmas}\label{sec:uni:lemma}
\begin{lem}
\label{uni:lem00}
The coefficients $\hat{\boldsymbol{\theta}}$ defined in~(\ref{uni:est2}) satisfies $\hat\theta_k = \sum_i d_{i,k} y_i$ with $d_{i,k} = o(1)$,  $1\leq k\leq c$.
\end{lem}
 {\it Proof of Lemma \ref{uni:lem00}:} It suffices to show every element of the matrix $(\mathbf{B}^T\mathbf{B}+\lambda^{\ast} \mathbf{D}^T\mathbf{D})^{-1}\mathbf{B}^T$ is $o(1)$. Because every column of $\mathbf{B}^T$ contains at most $p+1$ non-zero elements that sum to 1 by Lemma~\ref{uni:lem1}, it suffices to show that every element of the matrix $M^{-1}\BO{\Lambda}^{-1}=(\mathbf{B}^T\mathbf{B}+\lambda^{\ast} \mathbf{D}^T\mathbf{D})^{-1}$ is $o(1)$. Since $\BO{\Lambda}^{-1}$ is positive-definite, it suffices to show the diagonal elements of $M^{-1}\BO{\Lambda}^{-1}$ are $o(1)$.  For $1\leq i\leq c$, the largest eigenvalue  of $M^{-1}\BO{\Lambda}^{-1}$ is smaller than the largest eigenvalue of $(\mathbf{B}^T\mathbf{B})^{-1}$ since $ \mathbf{D}^T\mathbf{D}$ is positive semi-definite. By Lemma 2 in Zhou et al.\,(1998), the eigenvalues of $(\mathbf{B}^T\mathbf{B})^{-1}$ are $O(K/n)$. Hence the diagonal elements of $M^{-1}\BO{\Lambda}$ are all $O(K/n) = o(1)$.
\begin{lem}
\label{uni:lem1}
The B-splines satisfy $\sum_{k=1}^{K+p} B_k(x)=1$ for any $x\in(0,1)$.
\end{lem}
See page 201 in de Boor (1978).
\begin{lem} \label{uni:lem2}
 The B-splines with degree at least $1$ satisfy $\sum_{k=1}^{K+p} B_k(x)\{Kx-k+(p+1)/2\}=0 $ for any $x\in (0,1)$.
\end{lem}
{\it Proof of Lemma \ref{uni:lem2}:} By Lemma~\ref{uni:lem1}, $\sum_{k=1}^{K+p} B_k(x)\{Kx-k+(p+1)/2\}=0 $ is equivalent to
\begin{equation}
\sum_{k=1}^{K+p} B_k(x)k = Kx +(p+1)/2.
\label{uni:bspline1}
\end{equation}
We shall prove~(\ref{uni:bspline1}) by induction on $p$. Assume $p=1$. Let $k_x$ be the integer such that $x\in[k/K,(k+1)/K)$. Then $B_{k_x+1}(x) = -Kx+k+1$ and $B_{k_x+2}(x) = Kx-k$. It follows that
\begin{align*}
\sum_{k=1}^{K+1} B_k(x)k=&(-Kx+k_x+1)(k_x+1)+(Kx-k_x)(k_x+2)\\
=&(Kx -k_x)(k_x+2-k_x-1) +(k_x+1)\\
=&Kx + 1.
\end{align*}
Assume now the degree of the B-splines is $p$. We use $B_k^{[p]}(x)$ to denote the B-splines is of degree $p$. We use the recursive relation of de Boor,
\begin{align}
B_k^{[p]}(x) =& \frac{K}{p}\left[\left(x-\frac{k-p-1}{K}\right) B_{k-1}^{[p-1]}(x)+\left(\frac{k}{K}-x\right)B_{k}^{[p-1]}(x)\right]\nonumber\\
=&\frac{1}{p}\left[\left(Kx-k+p+1\right)B_{k-1}^{[p-1]}(x)-\left(Kx-k\right)B_{k}^{[p-1]}(x)\right].
\label{uni:deboor_recur}
\end{align}
It follows that
\begin{equation*}
\begin{split}
&p\left\{\sum_{k=1}^{K+p} B_k^{[p]}(x)k\right\}\\
=&\sum_{k=1}^{K+p}\left[\left(Kx-k+p+1\right)B_{k-1}^{[p-1]}(x)-\left(Kx-k\right)B_{k+1}^{[p-1]}(x)\right]k\\
=&\sum_{k=1}^{K+p-1}B_{k-1}^{[p-1]}(x)(Kx-k+p+1)k-\sum_{k=1}^{K+p-1}B_{k}^{[p-1]}(x)(Kx-k)k\\
=&\sum_{k=1}^{K+p-1}B_k^{[p-1]}(x)(Kx-k+p)(k+1)-\frac{1}{p}\sum_{k=1}^{K+p-1}B_k^{[p-1]}(x)(Kx-k)k\\
=&\sum_{k=1}^{K+p-1} B_k^{[p-1]}(x)(Kx-k+p+pk)\\
=&Kx + p+(p-1)\sum_{k=1}^{K+p-1}B_k^{[p-1]}(x)k \\
=&\left\{Kx + p +(p-1)(Kx+p/2)\right\}\\
=&p\left\{Kx + (p+1)/2\right\},
\end{split}
\end{equation*}
which is~(\ref{uni:bspline1}). Therefore, Lemma~\ref{uni:lem2} is proved.

\begin{lem}\label{uni:lem11}
Let $M=n/K$ be an integer. Let $\{B_1(x),\dots, B_c(x)\}$, where $c=K+p$,  be the the $B$-splines basis with knots $\{-p/K,-(p-1)/K,\dots, 0/K,1/K,\dots,K/K\}$. Then for $ k = q+1,\dots, K$,
 \begin{equation*}
 \sum_{i=1}^n B_{k}(x_i) = M
 \end{equation*}
 \end{lem}
  {\it Proof of Lemma \ref{uni:lem11}:} Proof by induction on $p$. Consider $p=0$. $B_{k}(x)=1$ if $x\in[k/K,(k+1)/K)$ and is 0 otherwise. So for fixed $k$, $B_{k}(x_i)=1$ if and only if $(i-1/2)/n\in[k/K, (k+1)/K)$, i.e., if and only if $i = nk/K+1, nk/K+1,\dots, n(k+1)/K$. Hence the case $p=0$ is proved. Now consider $p\geq 1$. By the recursive relation of de Boor in~(\ref{uni:deboor_recur}),
 \begin{align*}
  \sum_i B^{[p]}_{k}(x_i) =&\sum_i \frac{1}{p}\left[(Kx_i-k+p+1)B_{k-1}^{[p-1]}(x_r)-(Kx_i-k)B_{k}^{[p-1]}(x_i)\right]\\
  =&\frac{M(-k+p+1+k)}{p}+\frac{K}{p}\sum_i x_i\left\{B_{k-1}^{[p-1]}(x_i)-B_{k}^{[p]}(x_i) \right\}\\
  =&\frac{M(p+1)}{p} + \frac{K}{p}\left\{\sum_{i=1}^n x_i B_{k-1}^{[p-1]}(x_i)-\sum_{i=1}^{n-M} (x_i+1/K) B_{k-1}^{[p-1]}(x_i)\right\}\\
  =&\frac{M(p+1)}{p} + \frac{K}{p}\left\{\sum_{i=1}^n x_i B_{k-1}^{[p-1]}(x_i)-\sum_{r=1}^{n} (x_i+1/K) B_{k-1}^{[p-1]}(x_i)\right\}\\
  &\frac{M(p+1)}{p} + \frac{1}{p}\sum_{i=1}^n B_{k-1}^{[p-1]}(x_i)\\
  =&\frac{M(p+1)}{p} -\frac{1}{p}M\\
  =&M.
 \end{align*}
 So Lemma~\ref{uni:lem11} is proved. 

\begin{lem}\label{uni:lem0}
$P(1)=1, P^{\prime}(1)=p$.
\end{lem}
{\it Proof of Lemma~\ref{uni:lem0}:}
The expression of $P(x)$ in (\ref{uni:P3}) is rewritten here,
\begin{equation*}
P(x)= u_p +u_{p-1}x +\cdots +u_0 x^p +u_1 x^{p+1} +\cdots +u_p x^{2p}.
\end{equation*}
Hence,  $P(1)=2\sum_{i=1}^p u_i +u_0$ and $P^{\prime}(1)= p(2\sum_{i=1}^p u_i +u_0)$, so we only need to show that $2\sum_{i=1}^p u_i +u_0=1$. Let $\BO{C} =\BO{B}^T\BO{B}/M$. By  (\ref{uni:P4}), if $p<i<c-p$, then the coefficient vector $(u_p, u_{p-1},\cdots,u_0, u_1,\cdots, u_p)^T$ equals
 $(C_{i,i-p}, C_{i,i-p+1},\cdots,C_{i,i}, C_{i,i+1},\cdots, C_{i,i+p})^T$. Thus, $2\sum_{i=1}^p u_i +u_0=\sum_{|i-j|\leq p} C_{i,j}= \sum_j C_{i,j}$ because $C_{i,j}=0$ if $|i-j|>p$. Since $C_{i,j}= \sum_r B_i(x_r) B_j(x_r)/M$, $2\sum_{i=1}^p u_i +u_0= \sum_r \{B_i(x_r)\sum_j B_j(x_r)\}/M = \sum_r B_i(x_r)/M= 1$, where the last equality holds  by Lemma~\ref{uni:lem11}. 
\begin{lem}\label{uni:lemroot}
If $\{\psi_{1},\dots, \psi_m\}$ are the $m$ roots of $x^{2m}+(-1)^m=0$ satisfying that the real part of $\psi_{\nu}$ is positive,  then
\begin{equation}
\label{uni:rooteqn}
\prod_{j\neq \nu}(\psi_{\nu}^2-\psi_j^2)=(-1)^{m+1}m \psi_{\nu}^{-2}.
\end{equation}
\end{lem}
{\it Proof of Lemma~\ref{uni:lemroot}:} It is easy to see that $\{\psi_1^2,\dots, \psi_m^2\}$ are the $m$ roots of $x^m + (-1)^m=0$.
Thus, $\prod_{j=1}^m (x-\psi_j^2)=(-1)^m$. Taking derivative of $\prod_{j=1}^m (x-\psi_j^2)$ with respect to $x$ and letting  $x=\psi_{\nu}^2$ give~(\ref{uni:rooteqn}).
 \begin{lem}\label{uni:lem3}
Suppose $g(x) = \exp(- b |x|)$ with $b\neq 0$. 
\begin{equation*}
\begin{split}
\sum_{k,r} B_k(x)B_r(x_i)g(\frac{\bar{x}_k-\bar{x}_r}{h_n})= \left\{1-\frac{b}{Kh_n}\tilde{g}(x,x_i)+O\{(Kh_n)^{-2}\}\right\}g(\frac{x-x_i}{h_n}),
\end{split}
\end{equation*}
where 
\begin{align}
\label{uni:g_tilde}
\tilde{g}(x,x_i) =\begin{cases}
  2\sum_{k<r} B_k(x)B_r(x_i)(r-k) &\text{if} \quad x\geq x_i,\\
  2\sum_{k>r} B_k(x)B_r(x_i)(k-r)&\text{if} \quad x< x_i.
  \end{cases}
\end{align}
\end{lem}
{\it Proof of Lemma \ref{uni:lem3}:} Suppose  that $x\geq x_i$.
Take a Taylor expansion of $g(x)$ at the point $\frac{x-x_i}{h_n}$,
\begin{align*}
g(\frac{\bar{x}_k-\bar{x}_r}{h_n}) &= g(\frac{x-x_i}{h_n})\left\{1-\frac{b}{h_n} (|\bar{x}_k-\bar{x}_r|-|x-x_i|)+O\{(Kh_n)^{-2}\}\right\}\\
&=g(\frac{x-x_i}{h_n})\left\{1-\frac{b}{Kh_n} (|k-r|- Kx+Kx_i)+O\{(Kh_n)^{-2}\}\right\}.\\
\end{align*}
Hence if we drop the term $g(\frac{x-x_i}{h_n}) O\{(Kh_n)^{-2}\}$ in the above equality,

\begin{align*}
&\sum_{k,r} B_k(x)B_r(x_i)g(\frac{\bar{x}_k-\bar{x}_r}{h_n})\\
&=g(\frac{x-x_i}{h_n})\sum_{k,r} B_k(x)B_r(x_i)\left\{1-\frac{b}{Kh_n} (|k-r|- Kx+Kx_i)\right\}\\
&=g(\frac{x-x_i}{h_n})\left\{1-\frac{b}{Kh_n}\sum_{k,r}B_k(x)B_r(x_i)( |k-r|-Kx+Kx_i)\right\}\\
&=g(\frac{x-x_i}{h_n})\left\{1-\frac{b}{Kh_n}\sum_{k,r}B_k(x)B_r(x_i)\left( |k-r|-k+\frac{p+1}{2}+Kx_i\right)\right\}\\
&=g(\frac{x-x_i}{h_n})\left\{1-\frac{b}{Kh_n}\sum_{k,r}B_k(x)B_r(x_i)( |k-r|+r-k)\right\}\\
&=g(\frac{x-x_i}{h_n})\left\{1-\frac{2b}{Kh_n}\sum_{k<r}B_k(x)B_r(x_i)(r-k)\right\}.\\
\end{align*}
Note that in the above derivation, we used Lemma~\ref{uni:lem1} and \ref{uni:lem2}. The other case when $x<x_i$ can be similarly proved.


\begin{lem}
\label{uni:lem4}
The function $\tilde{g}$ defined in (\ref{uni:g_tilde}) satisfies
\[
\tilde{g}(x,x_i) = 0 \quad \text{if} \,\, |x-x_i|\geq (p+1)/K.
\]
\end{lem}
{\it Proof of Lemma \ref{uni:lem4}:} Suppose $x\geq x_i$. When $x-x_i\geq (p+1)/K$ and $k<r$, either $B_k(x)$ or $B_r(x_i)$ will be 0. The other case can be similarly proved. 
 \begin{lem}
 \label{uni:lem5}
Suppose $g(x) = \exp(- b |x|)$ with $b\neq 0$. 
\begin{equation*}
\sum_r B_r(x_i)g(\frac{r}{Kh_n})= \left[1+O\{(Kh_n)^{-1}\}\right]g(\frac{x_i}{h_n}).
\end{equation*}
\end{lem}
{\it Proof of Lemma \ref{uni:lem5}:} Take a Taylor expansion of $g(x)$ at the point $\frac{x_i}{h_n}$, 
\begin{align*}
g(\frac{r}{Kh_n}) &= g(\frac{x_i}{h_n})\left\{1-\frac{b}{h_n} \left(\frac{r}{K}-x_i\right)+O\{(Kh_n)^{-1}\}\right\}\\
&=g(\frac{x_i}{h_n})\left\{1-\frac{b}{Kh_n} (r-Kx_i)+O\{(Kh_n)^{-1}\}\right\}.\\
\end{align*}
Hence if we drop the term $g(\frac{x_i}{h_n}) O\{(Kh_n)^{-1}\}$ in the above equality,

\begin{align*}
\sum_r B_r(x_i)g(\frac{r}{Kh_n})&=g(\frac{x_i}{h_n})\sum_r B_r(x_i)\left\{1-\frac{b}{Kh_n} (r-Kx_i)\right\}\\
&=g(\frac{x_i}{h_n})\left\{1-\frac{b}{Kh_n}\sum_{r}B_r(x_i)(r-Kx_i)\right\}\\
&=g(\frac{x_i}{h_n})\left\{1-\frac{b}{Kh_n}\sum_{r}B_r(x_i)\frac{p+1}{2}\right\}\\
&=g(\frac{x_i}{h_n})\left(1-\frac{p+1}{Kh_n}\right).
\end{align*} 
 \begin{lem}
 \label{uni:lem6}
Assume $\psi$ is a complex number and $|\psi|=1$. For any nonnegative integer $\ell$,
$$
\int x^{\ell} e^{-\psi x} \mathrm{d} x = -e^{-\psi x}\sum_{k=1}^{\ell+1} \frac{ \ell! x^{\ell-k+1}}{(\ell-k+1)! } \bar{\psi}^k,$$
where $\bar{\psi}$ is the conjugate of $\psi$.
\end{lem}
{\it Proof of Lemma \ref{uni:lem6}:} The results of indefinite integrals of $\int x^{\ell} e^{a x} \cos (bx)\mathrm{d}x$ and $\int x^{\ell} e^{a x} \sin (bx)\mathrm{d}x$ are given by results 3 and 4 on page 230 of Gradshteyn and Ryzhik (2007). 

 \begin{lem}
 \label{uni:lem7}
Assume $|\psi|=1$ with positive real part.  For any nonnegative integer $\ell$,
$$
\int_0^{\infty} x^{\ell} e^{-\psi x} \mathrm{d}x = \ell!\bar{\psi}^{\ell+1},
$$
where $\bar{\psi}$ is the conjugate of $\psi$.
\end{lem}
{\it Proof of Lemma \ref{uni:lem7}:}  See Lemma~\ref{uni:lem6}. 

\begin{lem}
\label{uni:lem8}
If $\ell$ is even and $2\leq \ell\leq 2m-2$,
$$\sum_{\nu=1}^m \psi_{\nu}^{\ell} = 0.$$
\end{lem}
{\it Proof of Lemma \ref{uni:lem8}:} Assume $\{z_1,z_2,\dots, z_{2m}\}$ are all the roots of the equation $x^{2m}+(-1)^m=0$. Since $\ell$ is even, we can show that $\sum_{\nu=1}^m \psi_{\nu}^{\ell} = 1/2\sum_{i=1}^{2m} z_i^{\ell}$ because if $a +b\, \imath$ is a root of $x^{2m}+(-1)^m=0$, then $\pm a \pm b\, \imath $ are also roots.  Assume $m$ is odd first. Let $\omega = e^{\imath \pi/m}$. Note that $\omega$ is a primitive root of $x^{2m}=1$, and we can organize $\{z_1,\dots,z_{2m}\}$ in such a way that $z_i = \omega^i$. It follows that 
\[
\sum_{i=1}^{2m} z_i^{\ell} = \sum_{i=1}^{2m} \omega^{\ell i} = \omega^{\ell}\frac{1-\omega^{2m \ell}}{1-\omega^{\ell}}=0.
\]
For the case $m$ is even, let $\omega_0=e^{\imath \pi/(2m)}$.  We can also write $z_i = \omega_0^{1+2i}$, then
\[
\sum_{i=1}^{2m} z_i^{\ell} = \sum_{i=1}^{2m} \omega_0^{\ell(1+2 i)} = \omega_0^{\ell}\frac{1-\omega_0^{4m \ell}}{1-\omega_0^{2\ell}}=0.
\]

\begin{lem}
\label{uni:lem9}
$$\int_{-\infty}^{\infty} x^{\ell} H_m(x) \mathrm{d}x =\begin{cases} 
\hfill 1&:\quad \ell = 0\\
\hfill 0&:\quad \ell \,\mbox{ is odd}\\
 \hfill 0 &: \quad \ell\,\mbox{ is even and}\,\,2\leq  \ell\leq 2m-2\\
(-1)^{m+1} (2m)!&:\quad \ell = 2m
\end{cases} $$
\end{lem}
{\it Proof of Lemma \ref{uni:lem9}:}  Since $H_m(x)$ is symmetric about 0, the result for odd $\ell$ is obvious. Assume $\ell$ is even.  By Lemma~\ref{uni:lem7}, 
\begin{align*}
\int_{-\infty}^{\infty} x^{\ell} H_m(x)\mathrm{d}x &=\frac{1}{m}\sum_{\nu=1}^m \psi_{\nu}  \int_0^{\infty} x^{\ell} e^{-\psi_{\nu} x} \mathrm{d}x\\
&=\frac{\ell!}{m} \sum_{\nu=1}^m \psi_{\nu}\bar{\psi_{\nu}}^{\ell+1}\\
&=\frac{(-1)^{m+1}\ell!}{m}\sum_{\nu=1}^m \psi_{\nu}^{2m-\ell}.
\end{align*} 
If $\ell=0$, $\int_{-\infty}^{\infty}H_m(x)\mathrm{d}x = \frac{(-1)^{m+1}}{m}\sum_{\nu=1}^m \psi_{\nu}^{2m}=1$ as desired. If $\ell=2m$,  $\int_{-\infty}^{\infty} x^{2m} H_m(x)\mathrm{d}x =(-1)^{m+1} (2m)!$ also as desired. The case when $\ell$ is even and $2\leq \ell \leq 2m-2$ is proved by Lemma~\ref{uni:lem8}.

\section*{References}
\begin{description}

\item[]
\textsc{Claeskens G., Krivobokova T.}, and \textsc{Opsomer J.D.} (2009),
``Asymptotic properties of penalized spline estimators,"
\textit{Biometrika}, 96, 529-544.

\item[]
de \textsc{Boor, C.} (1978),
\textit{A Practical Guide to Splines},
Berlin: Springer.

\item[]
\textsc{Eilers, P.H.C.} and \textsc{Marx, B.D.} (1996),
``Flexbile smoothing with B-splines and penalties (with Discussion),"
\textit{Statist. Sci.}, 11, 89-121.

\item[]
\textsc{Gradshteyn I.S.} and \textsc{Ryzhik I.M.}(2007),
\textit{Table of Integrals, Series, and Products},
New York: Academic Press.

\item[]
\textsc{Hansen, B.E.} (2008),
``Uniform convergence rates for kernel estimation with dependent data,"
\textit{Econometric Theory}, 24, 726-748.

\item[]
\textsc{Holst, U., H\"{o}ssjer, O., Bj\"{o}rklund, C., Ragnarson, P.} and \textsc{Edner, H.} (1996),
``Locally weighted least squares kernel regression and statistical evaluation of LIDAR measurements,"
\textit{Environmetrics} 7: 401-416.

\item[]
\textsc{Kauermann, G., Krivobokova, T.} and \textsc{Fahrmeir, L.} (2009),
``Some asymptotic results on generalized penalized spline smoothing,''
\textit{J. R. Statist. Soc. Ser. B}, 71, 487-503.

\item[]
\textsc{Li,Y.} and \textsc{Ruppert D.} (2008),
``On the asymptotics of penalized splines,"
\textit{Biometrika}, 95, 415-436.

\item[]
\textsc{Messer, K.} and \textsc{Goldstein, L.}(1993),
``A new class of kernels for nonparametric curve estimation,''
\textit{Ann. Statist.}, 21, 179-195.
\item[]
\textsc{O'Sullivan, F.} (1986),
``A statistical perspective on ill-posed inverse problems (with discussion),''
\textit{Statist. Sci.}, 1, 505-527.

\item[]
\textsc{Opsomer, J.D.} and \textsc{Hall, P.} (2005),
``Theory for penalised spline regression,"
\textit{Biometrika}, 95, 417-436.

\item[]
\textsc{Ruppert, D., Wand, M.P.} and \textsc{Carroll, R.J.} (2003),
\textit{Semiparametric Regression},
 Cambridge: Cambridge University Press.

\item[]
\textsc{Ruppert, D., Wand, M.P., Holst, U.} and \textsc{ H\"{o}ssjer, O.} (1997),
``Local polynomial variance function estimation,"
\textit{Technometrics}, 39: 262-273.

\item[]
\textsc{Silverman, B.W.} (1984),
 ``Spline smoothing: the equivalent variable kernel method,''
 \textit{Ann. Statist.}, 12, 898-916.
 
 \item[]
\textsc{Stone, C.J.} (1980),
``Optimial rates of convergence for nonparametric estimators,''
\textit{Ann. Statist.}, 8, 1348-1360.

 \item[]
\textsc{Wand, M.P.} and \textsc{Jones, M.C.} (1995),
\textit{Kernel Smoothing},
London: Chapman \&Hall.

\item[]
\textsc{Wang X., Shen J.} and \textsc{Ruppert, D.} (2009),
``Local asymptotics of P-spline smoothing,'' 
\textit{EJS}, 4, 1-17.

\item[]
\textsc{Wang X.} and \textsc{Shen J.} (2010),
``A class of grouped Brunk estimators and penalized spline estimators for monotone regression,''
\textit{Biometrika}, 97, 585-601.

\item[]
\textsc{Zhou, S., Shen, X.} and \textsc{Wolfe, D.A.} (1998),
``Local asymptotics for regression splines and confidence regions,"
\textit{Ann. Statist.}, 26, 1760-1782.
\end{description}

\end{document}